\documentclass{birkjour}
  \usepackage{color}
  \usepackage[colorlinks=true,linkcolor=blue]{hyperref}
\def\cprime{$'$} \def\polhk#1{\setbox0=\hbox{#1}{\ooalign{\hidewidth
  \lower1.5ex\hbox{`}\hidewidth\crcr\unhbox0}}}
 \newcommand*{\al}{\alpha}
 \newcommand*{\ba}{\beta}
 \newcommand*{\Da}{\Delta}
 \newcommand*{\da}{\delta}
 \newcommand*{\Ga}{\Gamma}
 \newcommand*{\ga}{\gamma}
 \newcommand*{\ia}{\iota}
 \newcommand*{\ka}{\kappa}
 \newcommand*{\wk}{{\wt{\ka}}}
 \newcommand*{\kwk}{{\ka\wt{\ka}}}
 \newcommand*{\lda}{\lambda}
 \newcommand*{\na}{\nabla}
 \newcommand*{\om}{\omega}
 \newcommand*{\sa}{\sigma}
 \newcommand*{\ta}{\theta}
 \newcommand*{\ve}{\varepsilon}
 \newcommand*{\vp}{\varphi}
 \newcommand*{\BB}{{\mathbb B}}
 
 \newcommand*{\BE}{{\mathbb E}}
 \newcommand*{\BH}{{\mathbb H}}
 \newcommand*{\BJ}{{\mathbb J}}
 
 \newcommand*{\BM}{{\mathbb M}}
 \newcommand*{\BN}{{\mathbb N}}
 \newcommand*{\BR}{{\mathbb R}}
 \newcommand*{\BS}{{\mathbb S}}
 \newcommand*{\cA}{{\mathcal A}}
 \newcommand*{\cD}{{\mathcal D}}
 
 \newcommand*{\cH}{{\mathcal H}}
 \newcommand*{\cL}{{\mathcal L}}
 \newcommand*{\cR}{{\mathcal R}}
 \newcommand*{\cS}{{\mathcal S}}
 \newcommand*{\cU}{{\mathcal U}}
 \newcommand*{\cX}{{\mathcal X}}
 \newcommand*{\cZ}{{\mathcal Z}}
 \newcommand*{\gF}{{\mathfrak F}}
 \newcommand*{\gK}{{\mathfrak K}}
 \newcommand*{\gsm}{{\mathfrak m}}
 \newcommand*{\gN}{{\mathfrak N}}
 \newcommand*{\gsr}{{\mathfrak r}}
 \newcommand*{\gZ}{{\mathfrak Z}}
 \newcommand*{\sC}{{\mathsf C}}
 \newcommand*{\sP}{{\mathsf P}}
 \newcommand*{\sR}{{\mathsf R}}
 \newcommand*{\bal}{\begin{aligned}}
 \newcommand*{\eal}{\end{aligned}}
 \newcommand*{\uti}[1]{(#1)\space\space}
 \newcommand*{\mf}[1]{\boldsymbol{#1}}
 \newcommand*{\ci}{\mathaccent"7017 }   
 \newcommand*{\cqqncel}[2]{\ooalign{$\hfil#1\mkern1mu/\hfil$\crcr$#1#2$}}
 \newcommand*{\notsim}{\mathrel{\mathpalette\cqqncel\sim}}
 \newcommand*{\hb}[1]{$#1$}
 \newcommand*{\sdot}{\!\cdot\!}
 \newcommand*{\eea}[1]{\hbox{$[\![#1]\!]$}}
 \newcommand*{\vsdot}{\hbox{$\vert\sdot\vert$}}
 \newcommand*{\Vsdot}{\hbox{$\Vert\sdot\Vert$}}
 \newcommand*{\pe}{\hbox{$[{}\sdot{},{}\sdot{}]$}}
 \newcommand*{\pr}{\hbox{$(\cdot,\cdot)$}}
 \newcommand*{\prsn}{\hbox{$(\cdot\kern1pt|\kern1pt\cdot)$}}
 \newcommand*{\pw}{\hbox{$\dl{}\sdot{},{}\sdot{}\dr$}}
 \newcommand*{\mfa}{{\mf{a}}}
 \newcommand*{\mfW}{{\mf{W}}}
 \newcommand*{\mfrho}{{\mf{\rho}}}
 \newcommand*{\ol}{\overline}
 \newcommand*{\ul}{\underline}
 \newcommand*{\wh}{\widehat}
 \newcommand*{\wt}{\widetilde}
 \newcommand*{\btdot}{\mathrel{\raisebox{.25ex}{\hbox{\tiny$\scriptscriptstyle\bt$}}}}
 \newcommand*{\thBN}{\overset{{\hbox{\tiny$\scriptscriptstyle\bt$}}}{\BN}}
 \newcommand*{\thS}{\overset{{\hbox{\tiny$\scriptscriptstyle\bt$}}}{S}}
 
 \newcommand*{\vph}{\vphantom}
 \newcommand*{\hr}{\hookrightarrow}
 \newcommand*{\Llr}{\Longleftrightarrow}
 \newcommand*{\ra}{\rightarrow}
 \newcommand*{\dl}{\langle}
 \newcommand*{\dr}{\rangle}
 \newcommand*{\sdh}{\stackrel{d}{\hookrightarrow}}
 \newcommand*{\card}{\mathop{\rm card}\nolimits}
 \newcommand*{\cl}{{\rm cl}}
 \newcommand*{\diam}{\mathop{\rm diam}\nolimits}
 \newcommand*{\tdiv}{\mathop{\rm div}\nolimits}
 \newcommand*{\divgrad}{\mathop{\rm div\kern1pt grad}\nolimits}
 \newcommand*{\dom}{\mathop{\rm dom}\nolimits}
 \newcommand*{\grad}{\mathop{\rm grad}\nolimits}
 
 \newcommand*{\inj}{{\rm inj}}
 \newcommand*{\supp}{\mathop{\rm supp}\nolimits}
 \newcommand*{\Diff}{{\rm Diff}}
 \newcommand*{\Lis}{{\mathcal L}{\rm is}}
 \newcommand*{\loc}{{\rm loc}}
 \newcommand*{\is}{\subset}
 
 \newcommand*{\Chr}{\Ga_{ij}^k}
 \newcommand*{\whChr}{\wh{\Ga}_{ij}^k}    
 \newcommand*{\bt}{\bullet}
 \newcommand*{\es}{\emptyset}
 \newcommand*{\iy}{\infty}
 \newcommand*{\mt}{\mapsto}
 \newcommand*{\pl}{\partial}
 \newcommand*{\sh}{\sharp}
 \newcommand*{\Mg}{(M,g)}
 \newcommand*{\Mwhg}{(M,\wh{g})}
 \newcommand*{\Rmgm}{(\BR^m,g_m)}
 \newcommand*{\rc}{\hbox{r-c}\ }
 
 \newcommand*{\BUC}{BU\kern-.3ex C}
 \newtheorem{thm}{Theorem}[section]
 \newtheorem{cor}[thm]{Corollary}
 \newtheorem{lem}[thm]{Lemma}
 \newtheorem{prop}[thm]{Proposition}
 \theoremstyle{definition}
 \newtheorem{defn}[thm]{Definition}
 \newtheorem{rem}[thm]{Remark}
 \newtheorem{rems}[thm]{Remarks}
 \newtheorem{ex}[thm]{Example}
 \newtheorem{exs}[thm]{Examples}
 \theoremstyle{remark}
 \numberwithin{equation}{section}

\begin{document}
%
\firstpage{1} \volume{228} \Copyrightyear{2004} \DOI{003-0001}
%
%
%
%
%
%
%
\title[Function Spaces on Riemannian Manifolds]
 {Function Spaces on Uniformly Regular and Singular Riemannian Manifolds}
\author{Herbert Amann}
\address{Math. Institut\br
Universit\"at Z\"urich\br
Winterthurerstr. 190\br
Switzerland}
\email{herbert.amann@math.uzh.ch}
\date{\today}
\dedicatory{\large In memoriam Giuseppe Da Prato\\
who was a trailblazer of abstract evolution equations\\
and maximal regularity}

\begin{abstract}
This paper shows that the basic properties of Sobolev, Besov, and
Bessel potential spaces are valid on Riemannian manifolds with boundary,
which either have bounded geometry or posses singularities. In the latter
case the appropriate setting is that of Kondratiev-type weighted spaces.
The importance and usefulness of our results are indicated by a
demonstration
of a maximal regularity result for a linear parabolic initial value problem
on singular manifolds.
\end{abstract}

\maketitle
\section{Introduction}\label{sec-I}
It is well-known that Banach spaces of distributions---most notably Sobolev
and H\"older spaces---play a decisive role in the study of linear and
nonlinear differential equations. Wheras the theory of function spaces on
subdomains of Euclidean spaces is well developed, this is by far not true
if the underlying domain is a Riemannian manifold.

\smallskip
In recent years, the theory of differential equations on Riemannian manifolds
has found increasing interest. This is motivated both by intrinsic
differential geometric questions and by problems from applied fields like
fluid mechanics or numerical analysis, for example. In those connections
 evolution equations of parabolic type are of predominant importance.

\smallskip
During the last decades, the theory of linear and nonlinear parabolic
evolution equations on general Banach spaces has made great progress.
In particular, the local theory of quasilinear parabolic equations is by
now well established. It is based on linearization and so-called maximal
regularity results for linear equations (see, for example, \cite{Ama05a},
\cite[Chapter~5]{PS16a}).

\smallskip
In order that these abstract results become available for the study of
parabolic evolution equations in the global analysis setting, a~good
understanding of embedding, interpolation, point-wise multiplications,
and trace properties of Banach spaces of functions---more generally, of
sections of vector bundles---is fundamental. It is the purpose of this paper
to provide such results.

\smallskip
Basically, our paper consists of three parts. In the first one, which
comprises Sections \ref{sec-BG}--\ref{sec-F}, it is shown that, on
manifolds with boundary and bounded geometry, Sobolev, Besov, and Bessel
potential spaces are defined and possess the same properties as in the
well-known Euclidean setting.

\smallskip
In the second part we consider a class of singular manifolds. In this frame
weighted spaces of Kondratiev-type occur naturally. By an easy transposition
technique we prove that these spaces too possess all the properties
known to hold for the classical unweighted spaces. This is done in
Sections \ref{sec-S} and~\ref{sec-W}. Then, in Section~\ref{sec-D}, we
demonstrate, by a prototypical example, how the transposition method yields
maximal regularity results for uniformly parabolic equations on
singular manifolds.

\smallskip
It remains to produce concrete classes of singular manifolds.
This is achieved in the third part, in Sections
\ref{sec-M}--\ref{sec-K}.

\smallskip
In order to keep the exposition simple we present our results for function
spaces only, although everything applies to spaces of sections of general
tensor bundles. Also, as for applications, we consider merely singular
manifolds with smooth cuspidal point singularities. Cuspidal wedges and
corners are treated in~\cite{Ama25a} where an elaborate exposition of the
theory is provided. For this reason the proofs in the present paper are
rather brief and sometimes sketchy only.

\smallskip
So as to make this paper accessible to a broad readership, in the Appendix
we have collected the notations and conventions which we use throughout the
main body of this paper without further ado.
\section{Bounded Geometry}\label{sec-BG}
Let $\Mg$ be an \hbox{$m$-dimensional},
\ \hb{m\geq1}, Riemannian manifold with (possibly empty) boundary~$\pl M$.
For each
\hb{p\in\ci M:=M\setminus\pl M} and $v$ in $T_pM$ there exist a maximal open
interval~$J_p(v)$ about~$0$ in~$\BR$ and a unique geodesic
\hb{\ga_p(\cdot, v)\colon J_p(v)\ra M} satisfying
\hb{\ga_p(0,v)=p} and
\hb{\dot\ga_p(0,v)=v}, the maximal geodesic `starting at~$p$ in
direction~$v$'. The exponential map at~$p$ is defined by
\hb{\exp_p(v):=\ga_p(1,v)} for all
\hb{v\in T_pM} with
\hb{1\in J_p(v)}. Given any
\hb{p\in\ci M}, there exists
\hb{\rho(p)>0} such that $\exp_p$~is a diffeomorphism from the open ball
about the origin in~$T_pM$ with
radius~$\rho(p)$ onto an open neighborhood of~$p$ in~$\ci M$.
The supremum of all such~$\rho(p)$ is the injectivity radius,~$\inj(p)$,
of~$\ci M$ at~$p$. Given a nonempty subset~$S$ of~$\ci M$,
$$
\inj(S):=\inf_{p\in S}\inj(p)
$$
is the injectivity radius of~$S$. Note that $\inj(\ci M)$ may be zero.

\smallskip
Suppose $\pl M(\ve)$~is an open neighborhood of~$\pl M$ in~$M$ and $\chi$~is
a diffeomorphism from~$\pl M(\ve)$ onto
\hb{\pl M\times[0,\ve)} such that
\hb{\chi(p)=(p,0)} for
\hb{p\in\pl M}. Then $\bigl(\pl M(\ve),\chi\bigr)$ is a uniform collar
of~$\pl M$ in~$M$ of width~$\ve$. It is~a \emph{geodesic collar} if
\begin{equation}\label{BG.chi}
\chi^{-1}(p,t)=\ga_p\bigl(t,\nu(p)\bigr),
\ (p,t)\in\pl M\times[0,\ve),
\end{equation}
where $\nu$~is the inner (unit) normal of~$\pl M$.
Lastly, $\Mg$~has \emph{bounded curvature} if all covariant derivatives of
the Riemannian curvature tensor are bounded.
\begin{defn}\label{def-BG.B}
A~Riemannian manifold~$\Mg$ without boundary is said to have
\emph{bounded geometry} if it has a positive injectivity radius and bounded
curvature. If
\hb{\pl M\neq\es}, then $\Mg$~\emph{has bounded geometry} if
\begin{itemize}
\item[\rm(i)]
$\pl M$~has a uniform geodesic collar of width~$\ve$ in~$M$.
\item[\rm(ii)]
If
\hb{0<r<\ve}, then
\hb{\inj\bigl(M\setminus\pl M(r)\bigr)>0}.
\item[\rm(iii)]
$(\pl M,g_{\pl M})$~ has bounded geometry.
\item[\rm(iv)]
$\Mg$~has bounded curvature.
\end{itemize}
Here $g_{\pl M}$~is the restriction of~$g$ to the subbundle~$T\pl M$
of~$TM$.\qed
\end{defn}
Assume
\hb{\pl M\neq\es} and let $\na_{g_{\pl M}}$ be the Levi--Civita covariant
derivative of~$(\pl M,g_{\pl M})$. Th.~Schick~\cite{Schi01a} defines that
$\Mg$~has bounded geometry if
conditions (i), (ii), and (iv) of Definition~\ref{def-BG.B} apply,
$(\pl M,g_{\pl M})$ has a positive injectivity radius, and
all~$\na_{g_{\pl M}}$ covariant derivatives of the second fundamental form
are bounded. (Also see  B.~Ammann, N.~Gro{\ss}e, and V.~Nistor~\cite{AGN19c}
for a variant.)
\begin{thm}\label{thm-BG.BS}
Assume
\hb{\pl M\neq\es}. Then $\Mg$~has bounded geometry iff it has bounded
geometry in the sense of Schick.
\end{thm}
\begin{proof}
\cite[Theorem~XI.2.4.11]{Ama25a}.
\end{proof}
\section{Uniform Regularity}\label{sec-U}
To obtain flexible local descriptions of Riemannian manifolds with bounded
geometry we introduce the concept of uniformly regular Riemannian manifolds.

\smallskip
Let
\hb{U_\ka=\dom(\ka)} be the coordinate patch of a local chart~$\ka$ for~$M$.
Then $\ka$~is \emph{normalized}, provided
\begin{equation}\label{U.Q}
\ka(U_\ka)=Q_\ka^m
:=
\begin{cases}
(-1,1)^m,
&\text{if }U_\ka\is\ci M,\\
[0,1)\times(-1,1)^{m-1},
&\text{if }U_\ka\cap\pl M\neq\es.
\end{cases}
\end{equation}
An atlas~$\gK$ is \emph{normalized} if it consists of normalized charts.
A~normalized atlas is \emph{shrinkable} if there exists
\hb{r\in(0,1)} such that
\hb{\bigl\{\,\ka^{-1}(rQ_\ka^m)\ ;\ \ka\in\gK\,\bigr\}} is a covering of~$M$.
It has \emph{finite multiplicity} if there is
\hb{k\in\BN} such that any intersection of more than~$k$ coordinate patches
is empty.
\begin{defn}\label{def-U.K}
An atlas~$\gK$ for~$M$ is \emph{uniformly regular}~(ur) if
\begin{itemize}
\item[(i)]
it is normalized, shrinkable, and has finite multiplicity.
\item[(ii)]
\hb{\wk\circ\ka^{-1}\in\BUC^\iy\bigl(\ka(U_\kwk),\BR^m\bigr)} and\\
\hb{\|\wk\circ\ka^{-1}\|_{k,\iy}\leq c(k)},
\ \hb{\ka,\wk\in\gK},
\ \hb{k\in\BN}.\qed
\end{itemize}
\end{defn}
Here
\hb{U_\kwk:=U_\ka\cap U_\wk} is understood to be nonempty, and
\hb{\Vsdot_{k,\iy}} is the norm in~$\BUC^k$, the space of bounded and
uniformly continuous \hbox{$C^k$~functions}. Since $M$~is separable
and metrizable, it is not difficult to see that
a ur atlas is countable, that is, finite or countably infinite.
If $\gK$ and~$\wt{\gK}$ are atlases,
\hb{\gN(\ka,\wt{\gK}):=\{\,\wk\in\wt{\gK}\ ;\ U_\kwk\neq\es\,\}},
\ \hb{\ka\in\gK}.
\begin{defn}\label{def-U.KK}
Two \hbox{ur~atlases} $\gK$ and~$\wt{\gK}$ are \emph{equivalent} if
\begin{itemize}
\item[(i)]
\hb{\card\gN(\ka,\wt{\gK})+\card\gN(\wk,\gK)\leq c},
\ \hb{\ka\in\gK},
\ \hb{\wk\in\wt{\gK}}.
\item[(ii)]
\hb{\wk\circ\ka^{-1}\in\BUC^\iy\bigl(\ka(U_\kwk),\BR^m\bigr)} with\\
\hb{\|\wk\circ\ka^{-1}\|_{k,\iy}+\|\ka\circ\wk^{-1}\|_{k,\iy}\leq c(k)},
\ \hb{\ka\in\gK},
\ \hb{\wk\in\wt{\gK}},
\ \hb{k\in\BN}.\qed
\end{itemize}
\end{defn}
This induces an equivalence relation in the class of all ur atlases. Any
equivalence class is a \emph{ur structure} for~$M$.
A~\emph{ur manifold} is a manifold~$M$ together with a ur structure. Then
$\gK$~is a \emph{ur atlas for}~$M$ iff it belongs to this structure.
\begin{defn}\label{def-U.gg}
$\Mg$~is said to be a \emph{uniformly regular Riemannian (urR) manifold} if
there exists a ur atlas~$\gK$ for~$M$ such that
\begin{itemize}
\item[(i)]
\hb{\ka_*g\sim g_m,\ \ka\in\gK},
\item[(ii)]
\hb{\|\ka_*g\|_{k,\iy}\leq c(k)},
\ \hb{\ka\in\gK},
\ \hb{k\in\BN}.\qed
\end{itemize}
\end{defn}
Here $\ka_*g$~is the local representation of~$g$ on~$TQ_\ka^m$, the
push-forward by~$\ka$. Furthermore,
\hbox{${}\sim{}$}~holds uniformly w.r.t.\
\hb{\ka\in\gK}. It is immediate by Definition~\ref{def-U.KK} that this
determination is independent of the specific ur atlas~$\gK$.

\smallskip
We present a short list of easy examples on which we build below.
More sophisticated urR manifolds are found in Sections \ref{sec-M}--\ref{sec-K}.
\begin{exs}\label{exs-U.S}
(a)
$(\BR^m,g_m)$ and $(\BH^m,g_m)$, where
\hb{\BH^m:=\BR_+\times\BR^{m-1}}, are urR manifolds.

\smallskip
(b)
Compact manifolds are ur and all ur Riemannian metrics thereon are
equivalent.

\smallskip
(c)
Suppose that $(M_i,g_i)$,
\ \hb{i=1,2}, are ur and either $\pl M_1$ or~$\pl M_2$ is empty. Then
\hb{(M_1\times M_2,\,g_1+g_2)} is a urR manifold. Otherwise, it is a urR
manifold with corners.

\smallskip
(d)
Let
\hb{f\colon(M_1,g_1)\ra(M_2,g_2)} be an isometric diffeomorphism between
Riemannian manifolds. Then $(M_1,g_1)$ is ur iff $(M_2,g_2)$ is so. An
atlas~$\gK$ for~$M_1$ is ur iff
\hb{f_*\gK:=\{\,f_*\ka\ ;\ \ka\in\gK\,\}} is a ur atlas for~$M_2$.\qed
\end{exs}
The following fundamental result shows that urR manifolds yield local
descriptions of manifolds with bounded geometry.
\begin{thm}\label{thm-U.UB}
A Riemannian manifold is ur iff it has bounded geometry.
\end{thm}
The concept of urR manifolds has been introduced by the author
in~\cite{Ama12b}. In that paper it has been observed that a manifold without
boundary and bounded geometry is~ur. The converse is due to
M.~Disconzi, Y.~Shao, and G.~Simonett~\cite{DSS16a}. If
\hb{\pl M\neq\es}, then Theorem~\ref{thm-U.UB} is proved
in~\cite[Theorems XI.2.4.1 and XI.2.4.8]{Ama25a}.

\smallskip
Let $\Mg$ be ur and
\hb{\ve>0}. It is not difficult to see that there exists
a ur atlas~$\gK$ such that
\hb{\diam_g(U_\ka)<\ve} for
\hb{\ka\in\gK}.

\smallskip
An important technical property of ur manifolds is the subsequent fact. We
denote by~$\cD$ the space of smooth functions with compact support.
\begin{lem}\label{lem-U.LS}
Let $\gK$ be a ur atlas. There exist
\hb{\pi_\ka\in\cD\bigl(U_\ka,[0,1]\bigr)},
\ \hb{\ka\in\gK}, such that
\begin{itemize}
\item[(i)]
\hb{\{\,\pi_\ka^2\ ;\ \ka\in\gK\,\}} is a partition of unity subordinate
 to the open covering
 \hb{\{\,U_\ka\ ;\ \ka\in\gK\,\}} of~$M$.
\item[(ii)]
\hb{\|\ka_*\pi_\ka\|_{k,\iy}\leq c(k),\ \ka\in\gK,\ k\in\BN}.
\item[(iii)]
\hb{\chi\in\cD\bigl((-1,1)^m,[0,1]\bigr)} and
\hb{\chi|\supp(\ka_*\pi_\ka)=1} for
\hb{\ka\in\gK}.
\end{itemize}
\end{lem}
\begin{proof}
\cite[Lemma~3.2]{Ama12b}.
\end{proof}
The family
\hb{\bigl\{\,(\pi_\ka,\chi)\ ;\ \ka\in\gK\,\bigr\}} is said to be~a
\emph{localization system subordinate} to~$\gK$.
\section{Function Spaces}\label{sec-F}
In this section we introduce the most important Banach spaces of
distributions on urR manifolds and discuss their main properties.

\smallskip
Let $\Mg$ be a Riemannian manifold,
\hb{k\in\BN}, and
\hb{1\leq p\leq\iy}. Then
\begin{equation}\label{F.Wkp}
\|u\|_{k,p}=\|u\|_{k,p,M}
:=\sum_{j=0}^k\big\|\,|\na^ku|_{g_0^k}\,\big\|_{L_p(M)},
\ u\in C^k(M),
\end{equation}
where
\hb{\na=\na_g} denotes the Levi--Civita covariant derivative and
$L_p(M)$ equals $L_p(M,dV_g)$.

\smallskip
If
\hb{p=q<\iy}, then the \emph{Sobolev space}~$W_q^k(M)$,
\hb{1\leq q<\iy}, is the completion of $\cD(M)$
in~$L_{1,\loc}(M)$ with respect to the norm~%
\hb{\Vsdot_{k,q}}. The \emph{space $BC^k(M)$ of bounded
continuous
\hbox{$C^k$~functions}} on~$M$ is the Banach space of all
\hb{u\in C^k(M)} satisfying
\hb{\|u\|_{k,\iy}<\iy}, endowed with the norm~%
\hb{\Vsdot_{k,\iy}}. It follows that $W_q^0(M)$ equals $L_q(M)$ and
\hb{BC^0(M)=BC(M)}.

\smallskip
We write
\hb{\pr_{\ta,p}},
\ \hb{0<\ta<1},
\ \hb{1\leq p\leq\iy},
resp.~%
\hb{\pe_\ta},
\ \hb{0<\ta<1}, for the real, resp.\ complex, interpolation functors.

\smallskip
The \emph{Besov space}~$B_{p,r}^s(M)$ is defined for
\hb{s>0} and
\hb{1\leq p,r\leq\iy} by
\begin{equation}\label{F.Bqr}
B_{q,r}^s(M)
:=
\begin{cases}
\bigl(W_q^k(M),W_q^{k+1}(M)\bigr)_{s-k,\,r},
 &\text{if }k<s<k+1,\\
\bigl(W_q^k(M),W_q^{k+2}(M)\bigr)_{1/2,\,r},
 &\text{if }s=k+1,
\end{cases}
\end{equation}
if
\hb{p=q<\iy}, and by
$$
B_{\iy,r}^s(M)
:=
\begin{cases}
\bigl(BC^k(M),BC^{k+1}(M)\bigr)_{s-k,\,r},
 &\text{if }k<s<k+1,\cr
\bigl(BC^k(M),BC^{k+2}(M)\bigr)_{1/2,\,r},
 &\text{if }s=k+1,
\end{cases}
$$
where
\hb{k\in\BN}.

\smallskip
We introduce \emph{Bessel potential spaces}~$H_q^s(M)$ for
\hb{s\geq0} and
\hb{1<q<\iy} by
\begin{equation}\label{F.Hq}
H_q^s(M)
:=
\begin{cases}
\bigl[W_q^k(M),W_q^{k+1}(M)\bigr]_{s-k},
 &\text{if }k<s<k+1,\\
W_q^k(M),
 &\text{if }s=k.
\end{cases}
\end{equation}
As usual,
\hb{B_p^s(M):=B_{p,p}^s(M)},
\ \hb{1\leq p\leq\iy}. Then \emph{Slobodeckii spaces} are specified by
$$
W_q^s(M):=B_q^s(M),
\ s\in\BR_+\setminus\BN,
\ 1\leq q<\iy,
$$
and
\hb{\bigl[\,W_q^s(M)\ ;\ s\geq0\,\bigr]} is the
\emph{Sobolev--Slobodeckii space scale}. Also,
$$
BC^s(M):=B_\iy^s(M),
\ s\in\BR_+\setminus\BN,
$$
are the \emph{H\"older spaces}, and
\hb{\bigl[\,BC^s(M)\ ;\ s\geq0\,\bigr]} is the \emph{H\"older space scale}.
Note, however, that
\hb{B_\iy^k(M)\neq BC^k(M)} if
\hb{k\in\BN}.
\begin{rems}\label{rems-F.T}
(a)
For the sake of easy presentation we restrict ourselves to function spaces.
However, everything said in this paper applies to spaces of sections
of tensor bundles $(T_\tau^\sa M,g_\sa^\tau)$,
\ \hb{\sa,\tau\in\BN}. For this it suffices to replace \eqref{F.Wkp} by
\begin{equation}\label{F.WT}
\|u\|_{k,p,T_\tau^\sa M}
:=\sum_{j=0}^k\big\|\,|\na^ju|_{g_\sa^{\tau+j}}\,\big\|_{L_p(M)},
\ u\in C^k(T_\tau^\sa M).
\end{equation}
The preceding procedures then lead to the Sobolev
spaces~$W_q^k(T_\tau^\sa M)$, the spaces~$BC^k(T_\tau^\sa M)$, and,
consequently, to $B_{p,r}^s(T_\tau^\sa M)$ and~$H_q^s(T_\tau^\sa M)$.

\smallskip
(b)
It is clear from \eqref{F.WT} that
$$
\na\in\cL\bigl(\gF^{k+1}(T_\tau^\sa M),\gF^k(T_{\tau+1}^\sa M)\bigr),
\ \gF^k\in\{W_q^k,BC^k\}.
$$
Hence interpolation and \eqref{F.Bqr}, \eqref{F.Hq} yield
$$
\na\in\cL\bigl(\gF^{s+1}(T_\tau^\sa M),\gF^s(T_{\tau+1}^\sa M)\bigr)
$$
if $\gF^s$~belongs to
$$
\{\,W_q^s\ ;\ 1\leq q<\iy\,\}
\cup\{\,H_q^s\ ;\ 1<q<\iy\,\}
\cup\{\,B_{p,r}^s\ ;\ 1\leq p,r\leq\iy\,\}.\eqno{\qed}
 $$
\end{rems}
It should be noted that these definitions apply to any Riemannian manifold.
However, in such a generality these spaces are not too useful since they
may not possess important embedding and interpolation properties
(e.g.,~\cite{HeR08a}). The situation is different for urR manifolds.
This is due to the fundamental retraction-coretraction theorem below.
Thus we now suppose that
$$
\bt\quad
\Mg\text{ is a urR manifold}.
$$

\smallskip
Given a nonempty set~$S$,
\ \hb{\Ga(S):=\BR^S}. Let $\gK$ be a ur atlas. We introduce
$$
\BM_\ka:=
\begin{cases}
\BR^m,
&\text{if }U_\ka\is\ci M,\\
\BH^m,
&\text{otherwise},
\end{cases}
$$
and
$$
\mf{\Ga}(\BM):=\prod_{\ka\in\gK}\Ga(\BM_\ka).
$$

\smallskip
Suppose that
\hb{\bigl\{\,(\pi_\ka,\chi)\ ;\ \ka\in\gK\,\bigr\}} is a localization system
subordinate to~$\gK$. We set
$$
\gsr_\ka^cu:=\ka_*(\pi_\ka u),
\ u\in\Ga(M),
\quad \gsr_\ka v:=\pi_\ka\ka^*(\chi v),
\ v\in\Ga(\BM_\ka).
$$
Then
\begin{equation}\label{F.RC}
\cR^c\colon\Ga(M)\ra\mf{\Ga}(\BM),
\ u\mt(\gsr_\ka^cu)
\end{equation}
and
\begin{equation}\label{F.R}
\cR\colon\mf{\Ga}(\BM)\ra\Ga(M),
\ (v_\ka)\mt{\textstyle \sum_{\ka\in\gK}\gsr_\ka v_\ka}
\end{equation}
are linear mappings, and $\cR^c$~is a right inverse of~$\cR$.

\smallskip
If $\gF(\BM_\ka)$ is a linear subspace of~$\Ga(\BM_\ka)$
`of the same type for each~$\ka$', then
$$
\mf{\gF}(\BM):=\prod_\ka\gF(\BM_\ka)\is\mf{\Ga}(\BM).
$$
For example,
\hb{\mf{C}^k(\BM)=\prod_\ka C^k(\BM_\ka)}.

\smallskip
Let $X$ and~$Y$ be Banach spaces and
\hb{r\in\cL(X,Y)}. Then $r$~is a \emph{retraction} from~$X$ onto~$Y$ if it has
a continuous right inverse~$r^c$. In this case $(r,r^c)$~is called
\rc pair for~$(X,Y)$.

\smallskip
To have a unified presentation, we suppose that the
symbol~$\cX^s$ belongs to one of the following sets of symbols, where $\BUC$
means bounded and uniformly continuous.
$$
\bal
{}&\{\,\BUC^k\ ;\ k\in\BN\,\},
&&\{\,W_q^s\ ;\ 1\leq q<\iy,\ s\geq0\,\},\\
{}&\{\,H_q^s\ ;\ 1<q<\iy,\ s\geq0\,\},
&&\{\,B_{p,r}^s\ ;\ 1\leq p,r\leq\iy,\ s>0\,\},
\eal
$$
where
\hb{s=k} in the first case. The Banach spaces~$\cX^s(\BM_\ka)$ are particular
instances of the (much more general) spaces studied in great detail
in~\cite{Ama19a}.

\smallskip
We put, for
\hb{\mf{v}=(v_\ka)\in\mf{\cX}^s(\BM)},
$$
\|\mf{v}\|_{\ell_p(\mf{\cX}^s(\BM))}
:=\big\|\,(\|v_\ka\|_{\cX^s(\BM_\ka)})\,\big\|_{\ell_p},
\ 1\leq p\leq\iy.
$$
The linear subspace of~$\mf{\cX}^s(\BM)$,
$$
\ell_p\bigl(\mf{\cX}^s(\BM)\bigr)
:=\{\,\mf{v}\in\mf{\cX}^s(\BM)
\ ;\ \|\mf{v}\|_{\ell_p(\mf{\cX}^s(\BM))}<\iy\,\},
$$
is a Banach space with the norm~%
\hb{\Vsdot_{\ell_p(\mf{\cX}^s(\BM))}}.

\smallskip
The restriction of~$\cR$ to~$\ell_p\bigl(\mf{\cX}^s(\BM)\bigr)$,
resp.\ of~$\cR^c$ to~$\cX^s(M)$, where
$\cX^k(M)$ equals $BC^k(M)$, is again denoted by~$\cR$, resp.~$\cR^c$.
Then we can formulate the basic universal \rc theorem.
\begin{thm}\label{thm-F.RC}
Let $\Mg$ be ur. Then $(\cR,\cR^c)$ is an \rc pair for
\begin{itemize}
\item[(i)]
$\bigl(\ell_\iy(\mf{\BUC}^k(\BM)),BC^k(M)\bigr),
\ k\in\BN.$
\item[(ii)]
$\bigl(\ell_\iy(\mf{B}_{\iy,r}^s(\BM)),B_{\iy,r}^s(M)\bigr),
\ 1\leq r\leq\iy,\ s>0.$
\item[(iii)]
$\bigl(\ell_q(\mf{W}_q^s(\BM)),W_q^s(M)\bigr),
\ 1\leq q<\iy,\ s\geq0.$
\item[(iv)]
$\bigl(\ell_q(\mf{H}_q^s(\BM)),H_q^s(M)\bigr),
\ 1<q<\iy,\ s\geq0.$
\item[(v)]
$\bigl(\ell_q(\mf{B}_{q,r}^s(\BM)),B_{q,r}^s(M)\bigr),
 \ 1\leq q<\iy,\ 1\leq r\leq\iy,\ s>0.$
\end{itemize}
\end{thm}
\begin{proof}
Cases (i) and~(iii) with
\hb{s=k} are established by direct computations. The remaining statements are
then obtained by interpolation using the results
of \cite[Section~VI.2]{Ama19a}. (The reader may also
consult~\cite{Ama12b} where the basic ideas are already found.)
\end{proof}
\begin{rem}\label{rem-F.N}
Since $(\cR,\cR^c)$~is an \rc pair, it is well-known and easy to see that
\begin{equation}\label{F.LC}
u\mt\Bigl(\sum_\ka\|\gsr_\ka^cu\|_{W_q^s(\BM_\ka)}^q\Bigr)^{1/q}
\end{equation}
is a norm for~$W_q^s(M)$. Of course, it depends on the choice of the ur
atlas~$\gK$ and the localization system subordinate to it. However, since
$W_q^s(M)$~has been globally defined, another choice of a ur atlas and
a corresponding localization system yields an equivalent norm.
Analogous observations apply
to the spaces $BC^k(M)$, $H_q^s(M)$, and~$B_{p,r}^s(M)$.\qed
\end{rem}
One of the fundamental ramifications of the
 retraction-coretraction theorem
is the next statement.
\begin{thm}\label{thm-F.FM}
Let $\Mg$ be ur. The spaces $BC^k(M)$, $W_q^s(M)$, $H_q^s(M)$,
and $B_{p,r}^s(M)$ enjoy the same density, embedding,
and interpolation properties as the corresponding classical spaces
on~$(\BR^m,g_m)$, resp.\ $(\BH^m,g_m)$.
\end{thm}
\begin{proof}
This follows from the properties of \rc pairs, the features
of the associated function spaces on~$\BM_\ka$, and the properties of
Banach-space-valued sequence spaces.
\end{proof}
We content ourselves by presenting just one of the many consequences, namely,
Sobolev-type embedding theorems.
\begin{thm}\label{thm-F.S}
Assume that $\Mg$~is ur.
\begin{itemize}
\item[(i)]
{\rm (Sobolev)}
Suppose that
\hb{q_0,q_1\in[1,\iy)} and
\hb{0\leq s_0<s_1} satisfy
$$
s_1-m/q_1\geq s_0-m/q_0
$$
with a strict inequality unless
\hb{s_0,s_1\in\BN}. Then
\hb{W_{q_1}^{s_1}(M)\sdh W_{q_0}^{s_0}(M)}.
\item[(ii)]
{\rm (Morrey)}
If
\hb{0\leq t<s-m/q}, then
\hb{W_q^s(M)\hr BC^t(M)}.
\item[(iii)]
{\rm (Gagliardo--Nirenberg)}
Assume that
\hb{0\leq s_0<s<s_1},
\ \hb{1\leq q,q_0,q_1<\iy}, and
\hb{0<\ta<1}. Set
$$
s_\ta:=(1-\ta)s_0+\ta s_1,
\quad
1/q(\ta):=(1-\ta)/q_0+\ta/q_1.
$$
Let
$$
s-m/q=s_\ta-m/q(\ta),
\quad q\geq q(\ta),
$$
where
\hb{q_i>1} if
\hb{s_i\in\BN}. Then
\hb{W_{q_0}^{s_0}(M)\cap W_{q_1}^{s_1}(M)\hr W_q^s(M)} and
$$
\|u\|_{s,q}
\leq c\,\|u\|_{s_0,q_0}^{1-\ta}\,\|u\|_{s_1,q_1}^\ta.
$$
\end{itemize}
\end{thm}
\begin{proof}
\cite[Subsection~XII.3.3]{Ama25a}.
\end{proof}
With the help of the \rc pair $(\cR,\cR^c)$
the point-wise multiplier theorems of \cite[Section~VII.6]{Ama19a} can also
be lifted to the manifold~$\Mg$. We restrict ourselves to present the most
important of such results.

\smallskip
Let $\gF_i(M)$,
\ \hb{i=0,1,2}, be Banach spaces of regular distributions, that is,
\hb{\gF_i(M)\hr L_{1,\loc}(M)}. We write
$$
\gF_0(M)\btdot\gF_1(M)\hr\gF_2(M)
$$
if the point-wise product
\hb{(u_0,u_1)\mt u_0\btdot u_1} restricts to a continuous bilinear map from
\hb{\gF_0(M)\times\gF_1(M)} into~$\gF_2(M)$.
\begin{thm}\label{thm-F.FP}
Suppose that
\hb{s,t\geq0}. Then
\begin{equation}\label{F.MW}
BC^s(M)\btdot W_q^t(M)\hr W_q^t(M),
\ 1\leq q<\iy,
\end{equation}
provided either \hb{s>t} or
\hb{s=t\in\BN}, and
\begin{equation}\label{F.MB}
BC^s(M)\btdot BC^t(M)\hr BC^t(M),
\ 0\leq t\leq s.
\end{equation}
If
\hb{s>m/q}, then $W_q^s(M)$~is a multiplication algebra.
\end{thm}
Assume
\hb{\pl M\neq\es}. If $\Mg$~is ur, then so is $(\pl M,g_{\pl M})$.
Hence the spaces $BC^k(\pl M)$, $W_q^s(\pl M)$, $H_q^s(\pl M)$, and
$B_{p,r}^s(\pl M)$ are defined and the analogue of Theorem~\ref{thm-F.FM}
applies.

\smallskip
Let $\bigl(\pl M(\ve),\chi\bigr)$ be a geodesic collar of~$\pl M$. The normal
derivative of order
\hb{j\in\BN} of
\hb{u\in C^\iy(M)} at
\hb{p\in\pl M} is defined by, see~\eqref{BG.chi},
$$
\pl_\nu^ju(p)=\frac{\pl^ju}{\pl\nu^j}(p)
:=\Bigl(\frac d{dt}\Bigr)^j(u\circ\ga_p)(p,t)\Big|_{_{t=0}}.
$$
Thus
\hb{\pl_\nu^0u=u|\pl M}, the \emph{trace} of~$u$ on~$\pl M$.

\smallskip
The following \emph{trace theorem} is of predominant importance in the theory
of boundary value problems. It also explains why Besov spaces,
Slobodeckii spaces in particular, are of outstanding significance.
\begin{thm}\label{thm-F.T}
Suppose that $\Mg$~is ur and
\hb{\pl M\neq\es}. Also suppose that
$$
1<q<\iy,
\ j\in\BN,
\ s>j+1/q.
$$
Then
\hb{\vec\pl_\nu^j:=(\pl_\nu^0,\ldots,\pl_\nu^j)} is a retraction
from~$W_q^s(M)$ and from~$H_q^s(M)$ onto
$$
\prod_{i=0}^jB_q^{s-i-1/q}(\pl M).
$$
\end{thm}
\begin{proof}
Let $\cX^t$ be one of the symbols $W_q^t$ or~$H_q^t$. Denote by
$(\dot\cR,\dot\cR^c)$ the restriction of~$(\cR,\cR^c)$ to~$\pl M$
(in the obvious sense). It is seen that $(\dot\cR,\dot\cR^c)$ is an
\rc pair for
$$
\bigl(\ell_q(\mf{\cX}^t(\BR^{m-1})),\cX^t(\pl M)\bigr).
$$
Using it, we carry the (isotropic) half-space trace
Theorem~VIII.1.3.2 of~\cite{Ama19a} over to~$\pl M$.
\end{proof}
\begin{rems}\label{rems-F.TS}
(a)
Suppose we consider the more general case of tensor-valued sections
(Remark~\ref{rems-F.T}(a)). Then the spaces~$\cX^s(\BM_\ka)$ in
Theorem~\ref{thm-F.RC} have to be replaced by
$\cX^s(\BM_\ka,T_\tau^\sa\BR^m)$. Also, \eqref{F.MW} in the point-wise
multiplier theorem now reads
$$
BC^s(T_0^jM)\btdot W_q^t(T_{\tau+j}^\sa M)\hr W_q^t(T_\tau^\sa M)
$$
for
\hb{\sa,\tau,j\in\BN}. An analogous statement holds for~\eqref{F.MB}.
With these modifications everything said above
applies to the spaces of \hbox{$T_\tau^\sa$-sections}.

\smallskip
In fact, in~\cite{Ama25a} we admit general metric vector bundles (possibly
with Banach space fibers) equipped with a ur metric connection.
We also mention M.~Kohr and V.~Nistor~\cite{KoN22a} for some results on
Sobolev spaces on general Riemannian manifolds.

\smallskip
(b)
If $\ci W_q^s(M)$ is the closure of~$\cD(\ci M)$ in~$W_q^s(M)$, then
\hb{\ci W_q^s(M)=W_q^s(M)} for
\hb{0\leq s<1/q}, and
$$
\ci W_q^s(M)=\ker(\vec\pl_q^j),
\ j+1/q<s<j+1+1/q,
\ j\in\BN.
$$

\smallskip
(c)
Sobolev spaces of negative order are introduced
(in the \hbox{$T_\tau^\sa$-setting}) by
$$
W_q^{-s}(T_\tau^\sa M)
:=\bigl(\ci W_{q'}^s(T_\sa^\tau M)\bigr)',
\ s>0,
$$
w.r.t.\ the
\hb{L_{q'}(T_\sa^\tau M)\times L_q(T_\tau^\sa M)}-duality pairing, where
\hb{1/q+1/q'=1} and
\hb{1<q<\iy}. Then $(\cR,\cR^c)$ has a unique extension to these
spaces. This easily leads to interpolation and embedding theorems
for negative order spaces.\qed
\end{rems}
Nearly all papers on distribution spaces over Riemannian
manifolds are concerned with manifolds without boundary. Very often,
compactness is additionally required.

\smallskip
Bessel potential spaces have first been introduced by
R.~Strichartz~\cite{Str83a} on complete Riemannian manifolds without
boundary as the fractional power spaces of the Laplace--Beltrami operator.
Large parts of the theory of Triebel--Lizorkin, Sobolev, Besov, and Bessel
potential spaces have been lifted by
H.~Triebel \cite{Tri86a},~\cite{Tri87a} from~$\BR^m$ to Riemannian manifolds
without boundary and bounded geometry. A~unified presentation is given in
\cite[Chapter~7]{Tri92a}. He uses geodesic coordinates to
define Triebel-Lizorkin spaces by means of the analog of~\eqref{F.LC}.
Then Besov spaces are dealt with by real interpolation. For Bessel
potential spaces, Strichartz' technique is employed. The restriction to
geodesic coordinates does not allow to introduce local coordinate norms on
Besov spaces (see the introduction in Subsection~7.3.1 in~\cite{Tri92a}).
N.~Gro{\ss}e and C.~Schneider~\cite{GrSch13a} partly remove this restriction
by employing atlases which are uniformly equivalent to geodesic ones. This
then leads them to localized norms for Bessel potential and Besov
spaces~$B_q^s(M)$,
\ \hb{1<q<\iy}.

\smallskip
In the case of the H\"older--Zygmund spaces~$B_\iy^s$,
\ \hb{s>0}, Triebel's approach leads to unnatural restrictions on~$s$ which
are partially rectified by introducing global H\"older--Zygmund norms based
on difference norms along geodesics. This global access is definitely
restricted to manifolds without boundaries.

\smallskip
In contrast to all those works, we use a top down approach. This means:
the spaces are globally defined and, then, it is shown that they can be
localized by the $(\cR,\cR^c)$ pair. This method has been introduced
in~\cite{Ama12b} in the general context of tensor bundle sections
on a class of singular Riemannian manifolds which comprises the family of
Riemannian manifolds with bounded geometry having possibly a nonempty
boundary. As is witnessed by the results of the present section, this
approach is very flexible and allows to transfer practically all results,
known to hold in the Euclidean setting, to urR manifolds. In addition,
it has the advantage that it automatically guarantees that different
choices of the localization ingredients lead to
equivalent norms.

\smallskip
Besides its usefulness for lifting the function space results
from~$\BR^m$ to~$M$, the localization technique is fundamental for
establishing basic solvability results for (partial) differential equations
on urR manifolds (for example \cite{Ama16a}, \cite{Ama17b},
\cite{Ama20a},~\cite{Ama21a}.
Also see~\cite{AHS94a} for the earliest implementation of this technique
in the simplest case~$(\BR^m,g_m)$.)

\smallskip
As for Sobolev embedding theorems on Riemannian manifolds:
virtually all research concerns the special case
$$
W_q^1(M)\hr L_{q/(m-q)}(M).
$$
This is due to its
differential-geometric implications. Particularly noteworthy are the writings
of
Th.~Aubin \cite{Aub76a}, \cite{Aub82a}, \cite{Aub98a},
Th.~Aubin and Y.Y.~Li~\cite{AuL99a},
E.~Hebey \cite{Heb99a},~\cite{Heb06a},
E.~Hebey and F.~Robert~\cite{HeR08a}. Both Aubin and Hebey also study
the problem of optimal embedding constants and consider manifolds of bounded
geometry. Supplementary references are
E.~Hebey and M.~Vaugon~\cite{HeV01a} and
O.~Druet~\cite{Dru02a},
O.~Druet and E.~Hebey~\cite{DrH02a}, and
H.~Cheik Ali~\cite{Che22a} as well as the papers cited in these
publications.
There are many Sobolev estimates on Riemannian manifolds possessing a special
geometric structure. A~prototypical reference is
L.~Saloff-Coste~\cite{Sal02a}.

\smallskip
Not much seems to be known about general Sobolev embeddings and
Gagliardo--Nirenberg inequalities on noncompact manifolds e.g., L.~Adriano
and Ch. Xia~\cite{AdX10a}, N.~Badr~\cite{Bad09a}. The
situation is different if manifolds with additional geometric structures
(for example, curvature bounds, symmetries,~etc.) are considered. For this we
refer to the differential geometric literature. Also~\cite{Sal02a}
might be of interest.
\section{Singular Manifolds}\label{sec-S}
We now introduce a large class of urR manifolds by a conformal change of the
metric~$g$ with a suitable singularity function. In this section we prove
two technical lemmas which are crucial for the following.
\begin{defn}\label{def-S.rho}
Let $\Mg$ be a Riemannian manifold. Assume:
\begin{itemize}
\item[(i)]
\hb{\rho\in C^\iy\bigl(M,(0,1]\bigr)}.
\item[(ii)]
\hb{\wh{g}:=g/\rho^2} and
$\Mwhg$ is a urR manifold.
\item[(iii)]
\hb{\big\|\rho^{k+1}\,|\na^kd(\log\rho)|_{g_0^{k+1}}\,\big\|_{BC(M)}
 \leq c(k)},
 \ \hb{k\in\BN}.
\end{itemize}
Then $\rho$~is said to be~a
\emph{singularity function} for~$\Mg$, and $(M,g,\rho)$ is~a
\emph{sin\-gu\-lar Riemannian manifold} of conformal type, more precisely,
of type~$\rho$.\qed
\end{defn}
This is an updated version of the concept of singular manifolds
first introduced in~\cite{Ama12b}. It follows from Corollary~\ref{cor-W.eq}
below that condition~(iii) means that $\rho$~is an `admissible weight'
for~$\wh{g}$ in the denomination of
B.~Ammann, N.~Gro{\ss}e, and V.~Nistor~\cite{AGNi19a}.
\begin{rems}\label{rems-S.SN}
(a)
Definition~\ref{def-S.rho} implies that $M$~is ur.

\smallskip
(b)
Let $\rho_1$ and~$\rho_2$ be singularity functions for~$\Mg$ such that
\hb{\rho_1\sim\rho_2}. Then
\hb{g/\rho_1\sim g/\rho_2}. Denoting by~$\eea\rho$ the equivalence class
for~%
\hbox{${}\sim{}$} containing the representative~$\rho$, it is more precise
to say that $(M,g,\rho)$ is~a \emph{singular manifold of type}~%
\hbox{$\eea\rho$}, and
\hbox{$\eea\rho$}~is the
\emph{singularity type of}~$M$.

\smallskip
(c)
Note that
\hb{\wh{g}\sim g} iff
\hb{\rho\sim\mf1}. If
\hb{\rho\notsim\mf1}, then there exists a~sequence~$(p_j)$ which leaves every
compact subset of~$M$ such that $\bigl(\rho(p_j)\bigr)$ converges
to~$0$. Thus $\rho$~captures, in some sense, the
singular behavior of~$\Mg$ `near infinity' (of~$\Mwhg$).

\smallskip
(d)
In local coordinates,
$$
\wh{g}=\rho^{-2}g_{ij}dx^i\otimes dx^j,
\quad \wh{g}^*=\rho^2g^{ij}\frac\pl{\pl x^i}\otimes\frac\pl{\pl x^j}
$$
and
\hb{\sqrt{\wh{g}}=\rho^{-m}\sqrt{g}}. Consequently,
\hb{\vsdot_{\wh{g}_\sa^\tau}=\rho^{\tau-\sa}\,\vsdot_{g_\sa^\tau}}.\qed
\end{rems}
From now on \emph{our interest focuses on the singular
manifold}~$(M,g,\rho)$, whereas $\Mwhg$~is viewed as~a
`uniform regularization' (or `desingularization') thereof. Since $\wh{g}$~is
obtained by a conformal
change of~$g$, there is an intimate relationship between $\Mg$ and~$\Mwhg$.
In particular, the distribution spaces on the urR manifold~$\Mwhg$, defined
by means of the Levi--Civita covariant derivative
\hb{\wh{\na}:=\na_{\wh{g}}}, can be expressed solely in terms of
\hb{\na=\na_g} and~$\rho$. In this section we prepare the necessary technical
details.

\smallskip
We set
\hb{\gF(\wh{M}):=\gF\Mwhg} for a given Banach
space~$\gF$ of functions. If
\hb{a\in C(T_\tau^\sa M)}, then (see Remark~\ref{rems-F.T}(a))
$$
\|a\|_{BC(T_\tau^\sa M)}
:=\big\|\,|a|_{g_\sa^\tau}\,\big\|_{BC(M)}.
$$
We write, for
\hb{i,j\in\BN},
$$
a\in BC^\iy(T_j^iM;\rho)
\Llr\|\rho^{k+j-i}\na^ka\|_{BC(T_{j+k}^iM)}\leq c(k),
\ k\in\BN.
$$
\begin{lem}\label{lem-S.S}
Suppose that
\hb{k\geq1}. There exists
\begin{equation}\label{S.SBC}
S\in BC^\iy(T_2^1M;\rho)
\end{equation}
such that
\begin{equation}\label{S.NN}
\wh{\na}u-\na u=S^ku:=\sum_{t=1}^k\sC_{2+t}^1(S\otimes u),
\ u\in C^\iy(T_k^0 M),
\ k\geq1.
\end{equation}
\end{lem}
\begin{proof}
(1)
Let $\gK$ be a ur atlas for~$M$. In local coordinates,
\hb{\ka=(x^1,\ldots,x^m)},
$$
\wh{\na}u-\na u=(\whChr-\Chr)u_kdx^i\otimes dx^j,
\quad u=u_kdx^k,
$$
with the Christoffel symbols
$$
\Chr=\frac12g^{k\ell}(\pl_ig_{j\ell}+\pl_jg_{i\ell}-\pl_\ell g_{ij}),
\ \pl_i:=\frac\pl{\pl x^i}.
$$
An easy computation shows that
$$
\whChr-\Chr
=-\bigl(\da_i^k\pl_j+\da_j^k\pl_i-g_{ij}g^{k\ell}\pl_\ell\bigr)\log\rho.
$$
Using
\hb{g^{kj}g_{ji}=\da_i^k}, we get
$$
\da_i^k\frac\pl{\pl x^k}\otimes dx^i=\sC_1^2(g^*\otimes g).
$$
Thus, writing
\hb{S=S\bigl(d(\log\rho)\bigr)} and
$$
S:=-2d(\log\rho)\otimes\sC_1^2(g^*\otimes g)
+g\otimes\bigl(g^\sh d(\log\rho)\bigr)\in C^\iy(T_2^1M),
$$
we find
$$
\wh{\na}u-\na u=S\btdot u,
\ u\in C^\iy(T^*M).
$$
From \hb{\na g=0} it follows that
\begin{equation}\label{S.SG}
\rho^{\ell+1}\na^\ell S\bigl(d(\log\rho)\bigr)
=S\bigl(\rho^{\ell+1}\na^\ell d(\log\rho)\bigr),
\ \ell\in\BN.
\end{equation}
Remark~\ref{rems-S.SN}(d) implies
\hb{g_{ij}g^{k\ell}=\wh{g}^{ij}\wh{g}_{k\ell}}.  Hence,
by Definition~\ref{def-S.rho}(ii),
\,\hb{g_{ij}g^{k\ell}\sim1} uniformly w.r.t.\
\hb{\ka\in\gK}. From this and~\eqref{S.SG} we obtain
$$
\|\rho^{\ell+1}\na^\ell S\|_{BC(T_{\ell+2}^1M)}
\leq c\,\|\rho^{\ell+1}\na^\ell d(\log\rho)\|_{BC(T_{\ell+1}^0M)},
\ \ell\in\BN.
$$
This and Definition~\ref{def-S.rho}(iii) prove that \eqref{S.SBC} holds and
that \eqref{S.NN} applies for
\hb{k=1} with
\hb{S^1:=S}.

\smallskip
(2)
The assertion for
\hb{k\geq2} now follows by a straightforward somewhat tedious calculation.
\end{proof}
By means of this lemma we can express~$\wh{\na}^k$ in terms of~$\na^k$,
and vice versa:
\begin{lem}\label{lem-S.ab}
Let
\hb{k\geq1}. There exist
\begin{equation}\label{S.ab}
a_i^k,b_i^k\in BC^\iy(T_k^iM;\rho),
\ 0\leq i\leq k-1,
\end{equation}
such that
\begin{equation}\label{S.Na}
\wh{\na}^ku=\na^ku+\sum_{i=0}^{k-1}a_i^k\btdot\na^iu
\end{equation}
and
\begin{equation}\label{S.Nb}
\na^ku=\wh{\na}^ku+\sum_{i=0}^{k-1}b_i^k\btdot\wh{\na}^iu
\end{equation}
for
\hb{u\in C^\iy(M)}.
\end{lem}
\begin{proof}
It is easy to see that there exist
\begin{equation}\label{S.akB}
a^k\in BC^\iy(T_{k+1}^kM;\rho)
\end{equation}
such that
$$
S^ku=a^k\btdot u,
\ u\in C^\iy(T_k^0M).
$$
We proceed by induction.

\smallskip
(1)
If
\hb{k=1}, then
\hb{\wh{\na}^1=\na^1=d}. Hence we set
\hb{a_0^1:=0}. Assume
\hb{k=2}. Lemma~\ref{lem-S.S} shows that \eqref{S.Na} holds with
\hb{a_1^2:=S^1}.

\smallskip
(2)
Suppose that the assertions pertaining to \eqref{S.Na}
have already been proved for
\hb{0\leq j\leq k}. Then, by Lemma~\ref{lem-S.S},
$$
\bal
\wh{\na}^{k+1}u
&=\wh{\na}\Bigl(\na^ku+\sum_{i=0}^{k-1}a_i^k\btdot\na^iu\Bigr)\cr
&=(\na+S^k)\Bigl(\na^ku+\sum_{i=0}^{k-1}a_i^k\btdot\na^iu\Bigr)\cr
&=\na^{k+1}u+\sum_{i=0}^{k-1}\bigl(\na(a_i^k\btdot\na^iu)
 +a^k\btdot a_i^k\btdot\na^iu\bigr)+a^k\btdot\na^ku\cr
&=\na^{k+1}u+\sum_{i=0}^{k-1}(\na a_i^k+a^k\btdot a_i^k)\btdot\na^iu
 +\sum_{i=0}^{k-1}a_i^k\btdot\na^{i+1}u+a^k\btdot\na^ku.
 \eal
$$
Set
\hb{a_{-1}^k:=0},
\begin{equation}\label{S.aa1}
a_i^{k+1}:=\na a_i^k+a^k\btdot a_i^k+a_{i-1}^k,
\ 0\leq i\leq k-1,
\end{equation}
and
\hb{a_k^{k+1}:=a^k}. Then
\hb{\wh{\na}^{k+1}u=\na^{k+1}u+\sum_{i=0}^ka_i^{k+1}\btdot\na^iu}.

\smallskip
(3)
We show, by induction on~$k$, that
\begin{equation}\label{S.aki}
a_i^k\in BC^\iy(T_k^iM;\rho),
\ 0\leq i\leq k-1.
\end{equation}
If
\hb{k=2}, then
\hb{a_0^1=S^1} and \eqref{S.SBC} yield the claim.

\smallskip
Assume that \eqref{S.aki} applies for some
\hb{k\geq2}. Then
$$
\rho^{\ell+k+1-i}\na^\ell(\na a_i^k)
=\rho^{\ell+1+k-i}\na^{\ell+1}a_i^k,
\ \ell\in\BN,
$$
and the induction hypothesis guarantee that
\begin{equation}\label{S.Naki}
\|\rho^{\ell+k+1-i}\na^\ell(\na a_i^k)\|_{BC(T_{k+\ell+1}^iM)}\leq c(\ell)
\end{equation}
for
\hb{0\leq i\leq k-1}. Using the product rule, we find
$$
\rho^{\ell+k+1-i}\na^\ell(a^k\btdot a_i^k)
=\sum_{j=0}^\ell\binom\ell j\rho^{j+1}\na^ja^k
 \btdot\rho^{\ell-j+k-i}\na^{\ell-j}a_i^k.
$$
From \eqref{S.akB} it follows
$$
\|\rho^{j+1}\na^ja^k\|_{BC(T_{k+1+j}^kM)}\leq c(\ell),
\ 0\leq j\leq\ell<\iy.
$$
Assumption~\eqref{S.aki} yields
$$
\|\rho^{\ell-j+k-i}\na^{\ell-j}a_i^k\|_{BC(T_{k+\ell-j}^iM)}\leq c(\ell),
\ 0\leq j\leq\ell,
\ 0\leq i\leq k-1.
$$
This implies
\begin{equation}\label{S.aa2}
\|\rho^{\ell+k+1-i}\na^\ell(a^k
\btdot a_i^k)\|_{BC(T_{k+\ell+1}^iM)}\leq c(\ell),
\ \ell\in\BN,
\ 0\leq i\leq k-1.
\end{equation}
Also,
\hb{\rho^{\ell+k+1-i}\na^\ell a_{i-1}^k
   =\rho^{\ell+k-(i-1)}\na^\ell a_{i-1}^k} provides
\begin{equation}\label{S.i1}
\|\rho^{\ell+k+1-i}\na^\ell a_{i-1}^k\|_{BC(T_{k+\ell}^{i-1}M)}\leq c(\ell),
\ \ell\in\BN,
\ 0\leq i\leq k-1.
\end{equation}
Thus we obtain from \eqref{S.aa1} and
\hbox{\eqref{S.Naki}--\eqref{S.i1}} that
\begin{equation}\label{S.ki}
a_i^{k+1}\in BC^\iy(T_{k+1}^iM;\rho),
\ 0\leq i\leq k-1.
\end{equation}
It follows from
\hb{a_k^{k+1}=a^k} and \eqref{S.akB} that \eqref{S.ki} holds also for
\hb{i=k}. Thus the induction is complete and
\eqref{S.aki} is established. Due to step~(1), this proves the first claim.

\smallskip
(4)
Now we look at
$$
\wh{\na}^ju=\na^ju+\sum_{i=0}^{j-1}a_i^j\btdot\na^iu,
\ 0\leq j\leq k,
$$
as a lower triangular system of linear equations in the
unknowns~$\na^iu$ for
\hb{0\leq i\leq k}. Solving it by forward substitution, we verify that
the claim pertaining to \eqref{S.Nb} is also true.
\end{proof}
\section{Weighted Spaces}\label{sec-W}
Let
\hb{k\in\BN} and
\hb{\lda\in\BR}. Given
\hb{1\leq q<\iy}, we put
\begin{equation}\label{W.W}
\|u\|_{W_q^{k,\lda}(M;\rho)}
:=\sum_{i=0}^k\big\|\rho^{-\lda+i-m/q}\,|\na^iu|_{g_0^i}\,\big\|_{L_q(M)}.
\end{equation}
Also
\begin{equation}\label{W.B}
\|u\|_{BC^{k,\lda}(M;\rho)}
:=\sum_{i=0}^k\big\|\rho^{-\lda+i}\,|\na^iu|_{g_0^i}\,\big\|_{BC(M)}.
\end{equation}
Then the \emph{weighted Sobolev space}
$$
W_q^{k,\lda}(M;\rho)\text{ is the completion of }
\bigl(\cD(M),\Vsdot_{W_q^{k,\lda}(M;\rho)}\bigr)
\text{ in }L_{1,\loc}(M,dV_g).
$$
The \emph{weighted space of bounded smooth functions}
$$
\bal
&BC^{k,\lda}(M;\rho)\text{ is the Banach space of all }u\in C^k(M)\cr
\noalign{\vskip-1\jot}
&\text{for which \eqref{W.B} is finite, endowed
 with the norm }\Vsdot_{BC^{k,\lda}(M;\rho)}.
\eal
$$

\smallskip
\emph{Weighted Slobodeckii spaces} are defined for
\hb{1\leq q<\iy} and
\hb{s\in\BR_+\setminus\BN} by
$$
W_q^{s,\lda}(M;\rho)
:=\bigl(W_q^{k,\lda}(M;\rho),W_q^{k+1,\,\lda}(M;\rho)\bigr)_{s-k,\,q},
\ k<s<k+1.
$$
If we replace
\hb{\pr_{s-k,\,q}} by the complex interpolation functor~%
\hb{\pe_{s-k}}, then we obtain the \emph{weighted Bessel potential
spaces}~$H_q^{s,\lda}(M;\rho)$,
where
\hb{H_q^{k,\lda}:=W_q^{k,\lda}} for
\hb{k\in\BN}. \emph{Weighted Besov spaces} are defined analogously
to~\eqref{F.Bqr}. This yields, in particular, weighted H\"older spaces.

\smallskip
We set
\hb{W_q^k(M;\rho):=W_q^{k,0}(M;\rho)} and
\hb{BC^k(M;\rho):=BC^{k,0}(M;\rho)}. Also
\hb{W_q^k(\wh{M}):=W_q^k\bigl((M,\wh{g})\bigr)}, etc., and
\hb{L_q^\lda(M;\rho):=W_q^{0,\lda}(M;\rho)}.
\begin{thm}\label{thm-W.eq}
If
\hb{k\in\BN}, then
$$
W_q^k(M;\rho)\doteq W_q^k(\wh{M}),
\ 1\leq q<\iy,
$$
and
$$
BC^k(M;\rho)\doteq BC^k(\wh{M}).
$$
\end{thm}
\begin{proof}
Let
\hb{0\leq j\leq k} and set
\hb{a_j^j:=1}. It follows from
Remark~\ref{rems-S.SN}(d) and \eqref{S.Na} that
$$
\bal
|\wh{\na}^ju|_{\wh{g}_0^j}
&\leq\sum_{i=0}^j|a_i^j\btdot\na^iu|_{\wh{g}_0^j}
 =\sum_{i=0}^j|\rho^ja_i^j\btdot\na^iu|_{g_0^j}\cr
&\leq\sum_{i=0}^j|\rho^{j-i}a_i^j|_{g_0^{j-i}}
 \,|\rho^i\na^iu|_{g_0^i}.
\eal
$$
From this and \eqref{S.ab} we get
$$
|\wh{\na}^ju|_{\wh{g}_0^j}
\leq c(j)\sum_{i=0}^j|\rho^i\na^iu|_{g_0^i}.
$$
Now we employ Remark~\ref{rems-S.SN}(d) once more to deduce that
$$
\big\|\,|\wh{\na}^ju|_{\wh{g}_0^j}\big\|_{L_p(\wh{M})}
\leq c(j)\sum_{i=0}^j\big\|\rho^{i-m/p}
\,|\na^iu|_{g_0^i}\big\|_{L_p(M)},
\ 1\leq p\leq\iy.
$$
By summing these inequalities from~$0$ to~$k$ and using
\eqref{W.W}, \eqref{W.B}, we find that
$$
\|u\|_{W_q^k(\wh{M})}\leq c(k)\,\|u\|_{W_q^k(M;\rho)},
\quad \|u\|_{BC^k(\wh{M})}\leq c(k)\,\|u\|_{BC^k(M;\rho)}.
$$

\smallskip
Similar reasoning, based on \eqref{S.Nb}, shows
that the norm of $W_q^k(\wh{M})$, resp.~$BC^k(\wh{M})$,
is stronger than the one of
$W_q^k(M;\rho)$, resp.\ $BC^k(M;\rho)$. The theorem is proved.
\end{proof}
\begin{cor}\label{cor-W.eq}
Definition~\ref{def-S.rho}{\rm (iii)} is equivalent to
\hb{d(\log\rho)\in BC^\iy(T^*\wh{M})}.
\end{cor}
In order to deal with the case where
\hb{\lda\neq0} we need the subsequent commutator estimate.
\begin{lem}\label{lem-W.L}
Let
\hb{k\in\BN} and
\hb{\lda\in\BR}. Then
$$
\sum_{i=0}^k
|\wh{\na}^i(\rho^\lda u)|_{\wh{g}_0^i}
\sim\sum_{i=0}^k\rho^\lda\,|\wh{\na}^iu|_{\wh{g}_0^i}.
$$
\end{lem}
\begin{proof}
(1)
First we note that
\hb{\log\rho^\lda=\lda\log\rho} and the preceding corollary imply that,
setting
\hb{\da:=\rho^\lda},
$$
a:=d(\log\da)\in BC^\iy(T^*\wh{M}).
$$
Furthermore,
\hb{d\da=\da a}. Hence
\begin{equation}\label{W.d1}
\wh{\na}(\da u)=\da\wh{\na}u+d\da\otimes u=\da(\wh{\na}u+a\otimes u).
\end{equation}

\smallskip
(2)
We claim that
\begin{equation}\label{W.dk}
\wh{\na}^k(\da u)
=\da\Bigl(\wh{\na}^ku+\sum_{i=0}^{k-1}a_i^k\otimes\wh{\na}^iu\Bigr)
=\da\wh{\na}^ku+\sum_{i=0}^{k-1}a_i^k\otimes\da\wh{\na}^iu,
\ k\in\BN,
\end{equation}
where
\begin{equation}\label{W.aki}
a_i^k\in BC^\iy(T_{k-i}^0\wh{M}),
\ 0\leq i\leq k.
\end{equation}
This follows by induction, similarly as in step~(2) of the proof
of Lemma~\ref{lem-S.ab}. It yields
$$
|\wh{\na}^k(\da u)|_{\wh{g}_0^k}
\leq c(k)\sum_{i=0}^k\da\,|\wh{\na}^iu|_{\wh{g}_0^i},
\ k\in\thBN.
$$

\smallskip
(3)
As in the proof of Lemma~\ref{lem-S.ab}, we look at
$$
\wh{\na}^j(\da u)
=\da\wh{\na}^ju+\sum_{i=0}^ja_i^j\otimes\da\wh{\na}^iu,
\ 0\leq j\leq k,
$$
as a lower triangular system of linear equations in the
unknowns~$\da\wh{\na}^iu$,
\hb{0\leq i\leq k}. Solving it by forward substitution, we find
\hb{b_i^k\in BC^\iy(T_{k-i}^0\wh{M})},
\ \hb{0\leq i\leq k-1}, satisfying
$$
\da\wh{\na}^ku
=\wh{\na}^k(\da u)+\sum_{i=0}^{k-1}b_i^k\otimes\wh{\na}^i(\da u).
$$
From this we get
$$
|\da\wh{\na}^ku|_{\wh{g}_0^k}
\leq c(k)\sum_{i=0}^k|\wh{\na}^i(\da u)|_{\wh{g}_0^i},
\ k\in\BN.
$$
The assertion follows.
\end{proof}
We denote by $\rho^\lda W_q^s(\wh{M})$ the image space of~$W_q^s(\wh{M})$
under the map
\hb{u\mt\rho^\lda u}. Thus the regular distribution
\hb{u\in L_{1,\loc}(\wh{M})} belongs to $\rho^\lda W_q^s(\wh{M})$ iff
\hb{\rho^{-\lda}u\in W_q^s(\wh{M})}. Similar
definitions apply to the other spaces under consideration.
It is clear that the spaces $\rho^\lda W_q^s(\wh{M})$, etc., are Banach
spaces.
\begin{thm}\label{thm-W.WW}
Suppose that
\hb{\lda\in\BR}. The map
\hb{u\mt\rho^\lda u} is an isomorphism from
\begin{itemize}
\item[(i)]
$BC^s(\wh{M})$ onto $BC^{s,\lda}(M;\rho)$,
\ \hb{s\geq0}.
\item[(ii)]
$W_q^s(\wh{M})$ onto $W_q^{s,\lda}(M;\rho)$,
\ \hb{1\leq q<\iy},
\ \hb{s\geq0}.
\item[(iii)]
$H_q^s(\wh{M})$ onto $H_q^{s,\lda}(M;\rho)$,
\ \hb{1<q<\iy},
\ \hb{s\geq0}.
\item[(iv)]
$B_{p,r}^s(\wh{M})$ onto $B_{p,r}^{s,\lda}(M;\rho)$,
\ \hb{1\leq p,r\leq\iy},
\ \hb{s\geq0}.
\end{itemize}
\end{thm}
\begin{proof}
Lemma~\ref{lem-W.L} implies
$$
\sum_{j=0}^k\big\|\,|\wh{\na}^j(\rho^{-\lda}u)|_{\wh{g}_0^j}
\big\|_{L_p(\wh{M})}
\sim\sum_{j=0}^k\big\|\rho^{-\lda}\,|\wh{\na}^ju|_{\wh{g}_0^j}
\big\|_{L_p(\wh{M})}
$$
for
\hb{1\leq p\leq\iy}. It is a consequence of Lemma~\ref{lem-S.ab} that
the second sum is equivalent to
$$
\sum_{j=0}^k\big\|\rho^{-\lda}\,|\na^ju|_{\wh{g}_0^j}
\big\|_{L_p(\wh{M})}.
$$
We infer from Remark~\ref{rems-S.SN}(d) that
$$
\big\|\rho^{-\lda}\,|\na^ju|_{\wh{g}_0^j}\big\|_{L_p(\wh{M})}
=\big\|\rho^{-\lda+j-m/p}\,|\na^iu|_{g_0^j}\big\|_{L_p(M)}.
$$
These considerations show that
$$
\|u\|_{\rho^\lda W_q^k(\wh{M})}
\sim\|u\|_{W_q^{k,\lda}(M;\rho)},
\ k\in\BN,
\ 1\leq q<\iy,
$$
and
$$
\Vsdot_{\rho^\lda BC^k(\wh{M})}
\sim\Vsdot_{BC^{k,\lda}(M;\rho)}.
$$
This yields assertions (i) and~(ii) if
\hb{s=k}. The remaining claims now follow by interpolation.
\end{proof}
\begin{cor}\label{cor-W.WW}
The weighted Sobolev--Slobodeckii, Bessel potential, and Besov spaces
enjoy the same density, embedding, and interpolation properties as their
classical (non-weighted) counterparts on~$(\BR^m,g_m)$,
resp.\ $(\BH^m,g_m)$.
\end{cor}
\begin{rem}\label{rem-W.Wr}
Suppose that
\hb{\lda_1\geq\lda_0}. Then
$$
W_q^{s,\lda_1}(M)\hr W_q^{s,\lda_0}(M),
\ s\geq0,
\ 1\leq q<\iy.
$$
Analogous embeddings apply to Bessel potential and Besov spaces.
\end{rem}
\begin{proof}
If
\hb{s=k\in\BN}, then
\hb{\rho^{-\lda_0}=\rho^{-\lda_1}\rho^{\lda_1-\lda_0}\leq\rho^{-\lda_1}}
and \eqref{W.W} yield the assertion. The general case follows
by interpolation, and the rest is now clear.
\end{proof}
Besides the interpolation results which are an outflow of the preceding
corollary, there
is also the following theorem for the case of different parameters. For
simplicity, we consider only weighted
Sobolev--Slobodeckii spaces. A~corresponding statement applies to
Bessel potential spaces (with
\hbox{$\pr_{\ta,q}$} replaced by~%
\hbox{$\pe_\ta$}).
\begin{thm}\label{thm-W.I}
Suppose that
\hb{s_0,s_1\in\BR_+} with
\hb{s_0\neq s_1},
\ \hb{\ta\in(0,1)}, and $\lda_0$,~$\lda_1$ belong to~$\BR$. Set
\hb{s_\ta:=(1-\ta)s_0+\ta s_1} and
\hb{\lda_\ta:=(1-\ta)\lda_0+\ta\lda_1}. Then
$$
\bigl(W_q^{s_0,\lda_0}(M;\rho),W_q^{s_1,\lda_1}(M;\rho)\bigr)_{\ta,q}
\doteq W_q^{s_\ta,\lda_\ta}(M;\rho),
$$
provided either
\hb{s_\ta\notin\BN} or
\hb{q=2}.
\end{thm}
\begin{proof}
Let $\gK$ be a ur atlas for~$M$. Set
\hb{\rho_\ka:=\ka_*\rho(0)=\rho\bigl(\ka^{-1}(0)\bigr)}. It is not
difficult to see that
\hb{\ka_*\rho\sim\rho_\ka} uniformly w.r.t.\
\hb{\ka\in\gK}. Using this
`discretization', one shows that $(\cR,\cR^c)$ is an \rc pair for
$\bigl(\ell_q(\mfrho^\lda\mfW_q^s(\BM),\rho^\lda W_q^s(\wh{M})\bigr)$, where
$$
\mfrho^\lda\mfW_q^s(\BM)
:=\prod_{\ka\in\gK}\rho_\ka^\lda W_q^s(\BM_\ka).
$$
Since
\hb{W_q^s=B_q^s} if either
\hb{s\notin\BN} or
\hb{q=2}, we get the assertion by applying Corollary~VI.2.3.3(i)
and the isotropic version of
Theorems VII.2.7.2(i), VII.2.7.4, and VII.2.8.3 in~\cite{Ama19a}.
\end{proof}
\begin{rem}\label{rem-W.T}
Let
\hb{\sa,\tau\in\BN} and define
$$
\|u\|_{W_q^{k,\lda}(T_\tau^\sa M;\rho)}
:=\sum_{j=0}^k\big\|\rho^{-\lda+\tau-\sa+i-m/q}\ |\na^iu|_{g_\sa^{\tau+i}}
\,\big\|_{L_q(T_\tau^\sa M)}.
$$
An analogous specification yields ~%
\hbox{$\Vsdot_{BC^{k,\lda}(T_\tau^\sa M;\rho)}$}.
We introduce weighted Sobolev spaces,~$W_q^{k,\lda}(T_\tau^\sa M;\rho)$,
etc., by mimicking the definitions for the
spaces $W_q^{k,\lda}(M;\rho)$, etc.
Then the foregoing theorems apply in this setting also.\qed
\end{rem}
Of great significance in the theory of differential equations on Riemannian
manifolds are the continuity properties of the covariant derivative in the
function spaces under consideration. We content ourselves with exhibiting the
Sobolev--Slobodeckii space case, making now use of the preceding remark.
\begin{thm}\label{thm-W.Na}
Let
\hb{\lda\in\BR},
\ \hb{1\leq q<\iy},
\ \hb{\sa,\tau\in\BN}, and
\hb{s\geq0}. Then
$$
\na\in\cL\bigl(W_q^{s+1,\,\lda}(T_\tau^\sa M;\rho),
W_q^{s,\lda}(T_{\tau+1}^\sa M;\rho)\bigr).
$$
\end{thm}
\begin{proof}
Let
\hb{k\in\BN} and
\hb{u\in W_q^{k+1,\,\lda}(T_\tau^\sa M;\rho)}. Then
$$
\bal
\|\na u\|_{W_q^{k,\lda}(T_{\tau+1}^\sa M;\rho)}
 =\sum_{i=0}^k
 \big\|\rho^{-\lda+\tau+1-\sa+i-m/q}\,|\na^i\na u|_{g_\sa^{\tau+1+i}}
 \,\big\|_{L_q(M)}\qquad\cr
{}=\sum_{j=1}^{k+1}
 \big\|\rho^{-\lda+\tau-\sa-m/q}\,|\na^ju|_{g_\sa^{\tau+j}}
 \,\big\|_{L_q(M)}
 \leq\|u\|_{W_q^{k+1,\,\lda}(T_\tau^\sa M;\rho)}.
\eal
$$
This proves the assertion if
\hb{s=k}. If
\hb{s\notin\BN}, then we conclude by interpolation.
\end{proof}
As for point-wise multiplications: the following theorem is an easy
consequence of Theorems \ref{thm-F.FP} and~\ref{thm-W.WW} and
Remark~\ref{rems-F.TS}(a).
\begin{thm}\label{thm-W.M}
Let
\hb{\lda_0,\lda_1\in\BR},
\hb{1\leq q<\iy}, and
\hb{\sa,\tau,j\in\BN}. Then
$$
BC^{s_0,\lda_0}(T_0^jM;\rho)\btdot W_q^{s_1,\lda_1}(T_{\tau+j}^\sa M;\rho)
\hr W_q^{s_1,\,\lda_0+\lda_1}(T_\tau^\sa M;\rho)
$$
if either
\hb{s_0>s_1\geq0} or
\hb{s_0=s_1\in\BN}. Also
$$
BC^{s_0,\lda_0}(T_0^jM;\rho)\btdot BC^{s_1,\lda_1}(T_{\tau+j}^\sa M;\rho)
\hr BC^{s_1,\,\lda_0+\lda_1}(T_\tau^\sa M;\rho)
$$
with
\hb{\ 0\leq s_1\leq s_0}. If
\hb{s>m/q}, then
$$
W_q^{s,\lda}(M;\rho)\btdot W_q^{s,\lda}(M;\rho)
\hr W_q^{s,\,2\lda}(M;\rho).
$$
Thus $W_q^{s,\lda}(M;\rho)$ is a multiplication algebra for
\hb{s>m/q} iff
\hb{\lda=0}.
\end{thm}
We refrain from formulating the weighted analog of the trace
Theorem~\ref{thm-F.T} and refer instead to \cite{Ama12b} and~\cite{Ama25a}.
It should be noted that the present definition~\eqref{W.W} differs from the
one we employed in our earlier papers by the factor~$\rho^{-m/q}$.
This change allows for simpler formulations of the point-wise multiplier
theorems, for example.
\section{Differential Operators}\label{sec-D}
The importance of Theorem~\ref{thm-W.WW} is manifest by its corollary.
It is also of fundamental importance in the study of differential equations
on singular manifolds. Namely, it allows to carry over to the singular
setting all existence and regularity theorems which can be derived in the less
intricate frame of urR manifolds. Although this would lead too far afield, we
explain the procedure in a simple setting and show just one consequence.

\smallskip
Suppose that
\hb{s\geq0},
\ \hb{\lda,\mu\in\BR}, and
\hb{\ol{s}\geq s} with
\hb{\ol{s}>s} if
\hb{s\notin\BN}. Given
\begin{equation}\label{D.aa}
\mf{a}=(a_0,\ldots,a_k)
\in\prod_{j=0}^kBC^{\ol{s},\,\mu-\lda}(T_0^jM;\rho),
\ a_k\neq0,
\end{equation}
set
$$
\cA=\cA(\mf{a},\na)
:=\sum_{j=0}^ka_j\btdot\na^j,
$$
where
\hb{(a_j\btdot\na^j)u:=a_j\btdot(\na^ju)}. Then
\begin{equation}\label{D.ALW}
\cA\in\cL\bigl(W_q^{s+k,\,\lda}(M;\rho),W_q^{s,\mu}(M;\rho)\bigr)
\end{equation}
by Theorems \ref{thm-W.Na} and~\ref{thm-W.M}. Theorem~\ref{thm-W.WW}
guarantees that
$$
\sP_s^\lda
:=(u\mt\rho^\lda u)
\in\Lis\bigl(W_q^s(\wh{M}),W_q^{s,\lda}(M;\rho)\bigr).
$$
Hence
$$
\wh{\cA}:=\sP_s^{-\mu}\circ\cA\circ\sP_{s+k}^\lda
\in\cL\bigl(W_q^{s+k}(\wh{M}),W_q^s(\wh{M})\bigr).
$$
By means of \eqref{S.Nb} and \eqref{W.dk} we derive that
$$
\wh{\cA}=\cA(\wh{a},\wh{\na})=\sum_{j=0}^k\wh{a}_j\btdot\wh{\na}
$$
with
$$
\wh{\mf{a}}=(\wh{a}_0,\ldots,\wh{a}_k)
\in\prod_{j=0}^kBC^{\ol{s}}(T_0^j\wh{M}).
$$
Moreover, the map
\hb{\mfa\mt\wh{\mfa}} is an isomorphism. Given
\hb{f\in W_q^{s,\mu}(M;\rho)},
$$
u\in W_q^{s+k,\,\lda}(M;\rho)
\text{ satisfies }
\cA(\mf{a},\na)u=f
$$
iff
$$
\wh{u}=\sP_{s+k}^{-\lda}u\in W_q^{s+k}(\wh{M})
\text{ complies with }
\cA(\wh{a},\wh{\na})\wh{u}=\sP_s^{-\mu}f.
$$
The dependence of~$\wh{\mfa}$ on~$\mfa$ is rather intricate. However,
it is important that the leading coefficients~$a_k$ of~$\cA$ and $\wh{a}_k$
of~$\wh{\cA}$ enjoy the transparent relationship
\hb{\wh{a}_k=\rho^{\lda-\mu}a_k}.
\begin{ex}\label{ex-D.LB}
We consider the Laplace--Beltrami operator
\hb{\cA:=\Da=\divgrad} on~$\Mg$. We have to express it in terms of~$\na$.
For this we use that the divergence of tensor fields is the linear map
$$
\tdiv\colon C^\iy(T_\tau^\sa M)\ra C^\iy(T_\tau^{\sa-1}M),
\ \sa\in\thBN,
\ \tau\in\BN,
$$
defined by
$$
\tdiv a:=\sC_{\tau+1}^\sa(\na a),
\ a\in C^\iy(T_\tau^\sa M).
$$
Hence
$$
\bal
\Da u
=\sC_1^1\bigl(\na(g^\sh du)\bigr)
&=\sC_1^1\bigl(\na\sC_1^2(g^*\otimes du)\bigr)\cr
&=\sC_1^1\sC_1^2(g^*\otimes\na du)
=g^*\btdot\na^2u,
\eal
$$
since $\na$~commutes with contractions and
\hb{\na g^*=0}. Observe that
\hb{|g^*|_{g_2^0}=m}. This and Remark~\ref{rem-W.T} imply that
\hb{a_2\in BC^{\iy,\,-2}(T_0^2M;\rho)}. Thus
\hb{\mu=\lda-2} and
$$
\cA\in\cL\bigl(W_q^{s+2,\,\lda}(M;\rho),W_q^{s,\,\lda-2}(M;\rho)\bigr),
\ s\geq0.
$$
Consequently,
$$
\wh{\cA}=\sP_s^{2-\lda}\circ\cA\circ\sP_{s+2}^\lda
=\wh{a}_2\btdot\wh{\na}^2+\wh{a}_1\btdot\wh{\na}+\wh{a}_0.
$$
(An explicit computation shows that the lower order coefficients
are not zero.) Using
\hb{\wh{a}_2=\rho^{\lda-\mu}a_2=\rho^2g^*=\wh{g}^*}, we see that, setting
\hb{\wh{\Da}:=\Da_{\wh{g}}},
$$
\wh{\cA}=\wh{\Da}+\wh{a}_1\btdot\wh{\na}+\wh{a}_0,
\ \wh{a}_i\in BC^\iy(T_0^i\wh{M}).\eqno{\qed}
$$
\end{ex}
Mapping properties for elliptic and parabolic differential operators on urR
manifolds have been investigated in~\cite{Ama16a},~\cite{Ama17b}. Based on
those results and the preceding example, we can easily derive corresponding
statements in the weighted space setting. We content ourselves by giving
just one of the many possibilities.

\smallskip
Assume that
\hb{\pl M=\es} and fix a nontrivial compact subinterval~$J$ of\/~$\BR_+$
containing~$0$. We look at the initial value problem for the
reaction-diffusion equation on~$M$:
\begin{equation}\label{D.RD}
\pl_tu-\tdiv(a\btdot\grad u)=f\text{ on }M\times J,
\ u|_{t=0}=u_0\text{ on }M.
\end{equation}
Choose
\hb{\lda\in\BR} and
\hb{1<q<\iy} and assume that
\begin{equation}\label{D.fu}
(f,u_0)
\in L_q\bigl(J,L_q^{\lda-2}(M;\rho)\bigr)
\times W_q^{2-2/q,\,\lda-2/q}(M;\rho).
\end{equation}
Furthermore,
\begin{equation}\label{D.a}
a\in BC^1(T_1^1M)
\end{equation}
and there exists a constant
\hb{\ul{\ve}>0} such that
$$
(\eta|a\eta)_g\geq\ul{\ve}\,|\eta|_g^2,
\ \eta\in\Ga(TM),
$$
or, equivalently,
\begin{equation}\label{D.ax}
(\xi|a\xi)_{g^*}\geq\ul{\ve}\,|\xi|_{g^*}^2,
\ \xi\in\Ga(T^*M).
\end{equation}
This means that
\begin{equation}\label{D.A}
\cA:=-\tdiv(a\btdot\grad)
\end{equation}
is a uniformly (strongly) elliptic differential operator on~$C^2(M)$.
\begin{thm}\label{thm-D.D}
Problem~\ref{D.RD} has for each~$(f,u_0)$ satisfying \eqref{D.fu}
a unique solution
\begin{equation}\label{D.u}
u=u(f,u_0)
\in L_q\bigl(J,W_q^{2,\lda}(M;\rho)\bigr)
\cap W_q^1\bigl(J,L_q^{\lda-2}(M;\rho)\bigr).
\end{equation}
The map
\hb{(f,u_0)\mt u(f,u_0)} is linear and continuous.
\end{thm}
\begin{proof}
(1)
Using local coordinates, for example, it is verified that
$$
-\cA
=(g^\sh a)\btdot\na^2+g^\sh\tdiv(a)\btdot\na
=a_2\btdot\na^2+a_1\btdot\na,
$$
where
\hb{g^\sh a=\sC_1^2(g^*\otimes a)\in C^1(T_0^2M)}. We find
$$
|a_2|_{g_2^0}=|g^\sh a|_{g_2^0}
\leq|g^*|_{g_2^0}\,|a|_{g_1^1}
=m\,|a|_{g_1^1}\in BC(M).
$$
This implies
$$
a_2\in BC^{0,\,-2}(T_0^2M;\rho).
$$
Similarly, by \eqref{D.a},
$$
|a_1|_g=|g^\sh\tdiv(a)|_g
\leq m\,|\tdiv a|_{g^*}\leq m\,|\na a|_{g_1^2}\in BC(M),
$$
which yields
$$
\rho^2(\rho^{-1}\,|a_1|_g)=\rho\,|a_1|_g\leq|a_1|_g\in BC(M).
$$
Hence
$$
a_1\in BC^{0,\,-2}(T_0^1M;\rho).
$$
As in Example~\ref{ex-D.LB}, the leading coefficient of
\hb{-\wh{\cA}=-\sP_0^{2-\lda}\cA\sP_2^\lda} reads
$$
\wh{a}_2=\rho^2a_2=\rho^2g^\sh a=\wh{g}^\sh a\in BC(T_0^2\wh{M}).
$$
It follows from~\eqref{D.ax} that
$$
\bal
\wh{a}_2\btdot(\xi\otimes\xi)
&=\dl\xi,\wh{a}_2\xi\dr
 =\rho^2\dl\xi,(g^\sh a)\xi\dr\cr
&=\rho^2(\xi|a\xi)_{g^*}
 \geq\ul{\ve}\rho^2\,|\xi|_{g^*}^2
 =\ul{\ve}\,|\xi|_{\wh{g}^*}^2
\eal
$$
for
\hb{\xi\in\Ga(T^*M)}. Hence $\wh{\cA}$~is a uniformly elliptic operator
on~$C^2(\wh{M})$.

\smallskip
(2)
For abbreviation,
$$
\wh{W}_q^s:=W_q^s(\wh{M}),
\quad W_q^{s,\lda}:=W^{s,\lda}(M;\rho).
$$
By Corollary~\ref{cor-W.WW} and Remark~\ref{rem-W.Wr},
\begin{equation}\label{D.WL}
W_q^{2,\lda}\sdh L_q^{\lda-2}.
\end{equation}
Theorems \ref{thm-F.FM} and~\ref{thm-W.I} imply
\begin{equation}\label{D.LW}
(\wh{L}_q,\wh{W}_q^2)_{1-1/q,\,q}\doteq\wh{W}_q^{2-2/q},
\quad (L_q^{\lda-2},W_q^{2,\lda})_{1-1/q,\,q}\doteq W_q^{2-2/q,\,\lda-2/q},
\end{equation}
respectively. Since, by Theorem~\ref{thm-W.WW},
\begin{equation}\label{D.PP}
\sP_0^{\lda-2}\in\Lis(\wh{L}_q,L_q^{\lda-2}),
\quad \sP_2^\lda\in\Lis(\wh{W}_q^2,W_q^{2,\lda}),
\end{equation}
we obtain from \eqref{D.LW} by interpolation, see Theorem~\ref{thm-W.I}, that
\begin{equation}\label{D.WLW}
\bal
&\sP_{2-2/q}^{\lda-2/q}\text{ is an isomorphism from }
\wh{W}_q^{2-2/q}\doteq(\wh{L}_q,\wh{W}_q^2)_{1-1/q,\,q}\cr
\noalign{\vskip-1\jot}
&\text{onto }
 (L_q^{\lda-2},W_q^{2,\lda})_{1-1/q,\,q}\doteq W_q^{2-2/q,\,\lda-2/q}.
\eal
\end{equation}
Hence assumption~\eqref{D.fu} yields
\begin{equation}\label{D.un}
\wh{u}_0:=\sP_{2-2/q}^{2/q-\lda}u_0\in(\wh{L}_q,\wh{W}_q^2)_{1-1/q,\,q}.
\end{equation}

\smallskip
(3)
The point-wise extension of
\hb{\sP_s^\lda\in\Lis(\wh{W}_q^s,W_q^{s,\lda})} over~$J$ is again denoted
by~$\sP_s^\lda$. Then we get from \eqref{D.PP} that
$$
\sP_0^{\lda-2}\in\Lis\bigl(L_q(J,\wh{L}_q),L_q(J,L_q^{\lda-2})\bigr)
$$
and
\begin{equation}\label{D.PJ}
\sP_2^\lda\in\Lis\bigl(L_q(J,\wh{W}_q^2),L_q(J,W_q^{2,\lda})\bigr).
\end{equation}
Thus \eqref{D.u} implies
\begin{equation}\label{D.ff}
\wh{f}:=\sP_0^{2-\lda}f\in L_q(J,\wh{L}_q).
\end{equation}
Step~(1) guarantees that
\hb{\pl_t+\wh{\cA}} is a uniformly parabolic differential operator
on~$C^2(\wh{M})$. Thus \eqref{D.un}, \eqref{D.ff}, and Theorem~1.23
in~\cite{Ama17b} guarantee that the initial value problem
\begin{equation}\label{D.I}
\pl_tv+\wh{\cA}v=\wh{f}\text{ on }M\times J,
\quad v|_{t=0}=\wh{u}_0\text{ on }M
\end{equation}
has a unique solution
\begin{equation}\label{D.ufu}
\wh{u}=\wh{u}(\wh{f},\wh{u}_0)
\in L_q(J,\wh{W}_q^2)\cap W_q^1(J,\wh{L}_q)
\end{equation}
and that
\hb{(\wh{f},\wh{u}_0)\mt\wh{u}} is linear and continuous.
Assertion~\eqref{D.ufu} is equivalent to
\begin{equation}\label{D.udu}
(\wh{u},\pl_t\wh{u})\in L_q(J,\wh{W}_q^2)\times L_q(J,\wh{L}_q).
\end{equation}

\smallskip
(4)
We set
\begin{equation}\label{D.uPu}
u:=\sP_2^\lda\wh{u}.
\end{equation}
Then, by \eqref{D.WL}, \eqref{D.PJ}, and~\eqref{D.udu},
\begin{equation}\label{D.u2}
u\in L_q(J,W_q^{2,\lda})\hr L_q(J,L_q^{\lda-2}).
\end{equation}
Moreover, see \eqref{D.ff},
$$
(\pl_t+\wh{\cA})\wh{u}
=(\pl_t+\sP_0^{2-\lda}\cA\sP_2^\lda)\wh{u}
=\sP_0^{2-\lda}f\in L_q(J,\wh{L}_q)
$$
and \eqref{D.uPu} imply
\begin{equation}\label{D.dPu}
\pl_t(\sP_2^{-\lda}u)
=\sP_2^{-\lda}\pl_tu
=\sP_0^{2-\lda}(-\cA u+f)\in L_q(J,\wh{L}_q).
\end{equation}
Note that
\hb{0<\rho\leq1} yields
$$
\bal
\|\pl_tu\|_{L_q(J,L_q^{\lda-2})}
&=\|\rho^{2-\lda-m/q}\pl_tu\|_{L_q(J,L_q)}\cr
&\leq\|\rho^{-\lda-m/q}\pl_tu\|_{L_q(J,L_q)}
 =\|\rho^{-\lda}\pl_tu\|_{L_q(J,\wh{L}_q)}.
\eal
$$
From this and \eqref{D.dPu} it follows that
\begin{equation}\label{D.dtu}
\bal
\|\pl_tu\|_{L_q(J,L_q^{\lda-2})}
&\leq \|\sP_0^{2-\lda}(-\cA u+f)\|_{L_q(J,\wh{L}_q)}\cr
&\leq c(\|\cA u\|_{L_q(J,L_q^{\lda-2})}+\|f\|_{L_q(J,L_q^{\lda-2})}\cr
&\leq c\bigl(\|u\|_{L_q(J,W_q^{2,\lda})}+\|f\|_{L_q(J,L_q^{\lda-2})}\bigr),
\eal
\end{equation}
where we also used \eqref{D.ALW} and step~(1). Consequently, we get from
\eqref{D.u2} that
$$
u\in L_q(J,W_q^{2,\lda})\cap W_q^1(J,L_q^{\lda-2}).
$$

\smallskip
(5)
We know from \cite[Theorem~III.4.10.2]{Ama95a} that
$$
L_q(J,\wh{W}_q^2)\cap W_q^1(J,\wh{L}_q)
\hr\BUC\bigl(J,(\wh{L}_q,\wh{W}_q^2)_{1-1/q,\,q}\bigr)
=:\wh{\BE}
$$
and
$$
L_q(J,\wh{W}_q^{2,\lda})\cap W_q^1(J,L_q^{\lda-2})
\hr\BUC\bigl(J,(L_q^{\lda-2},W_q^{2,\lda})_{1-1/q,\,q}\bigr)
=:\wh{\BE}^\lda.
$$
Hence, setting
\hb{\ga_0\wh{u}:=\wh{u}|_{t=0}},
$$
\ga_0\wh{u}=\lim_{t\ra0}\wh{u}(t)\text{ in }\wh{\BE}.
$$
We deduce from \eqref{D.WLW} that
$$
\sP_{2-2/q}^{\lda-2/q}\ga_0\wh{u}
=\lim_{t\ra0}\sP_{2-2/q}^{\lda-2/q}\wh{u}(t)
=\lim_{t\ra0}u(t)=\ga_0u\text{ in }\BE^\lda.
$$
Consequently,
$$
\ga_0u=\sP_{2-2/q}^{\lda-2/q}\ga_0\wh{u}=u_0,
$$
due to \eqref{D.un}.

\smallskip
Since the last part of the assertion follows from \eqref{D.uPu},
\eqref{D.dtu}, and the linearity and continuity of
\hb{(\wh{f},\wh{u}_0)\mt\wh{u}}, the theorem is proved.
\end{proof}
Let $E_0$ and~$E_1$ be Banach spaces with
\hb{E_1\sdh E_0}. Then $\cH(E_1,E_0)$ denotes the set of all
\hb{A\in\cL(E_1,E_0)} such that~$-A$, considered as a
linear operator in~$E_0$ with domain~$E_1$, is the infinitesimal generator of
a strongly continuous analytic semigroup
\hb{\{\,e^{-tA}\ ;\ t\geq0\,\}} on~$E_0$, that is, in~$\cL(E_0)$.
\begin{thm}\label{thm-D.A}
\hb{\cA\in\cH\bigl(W_q^{2,\lda}(M;\rho),L_q^{\lda-2}(M;\rho)\bigr)}.
\end{thm}
\begin{proof}
Theorem~\ref{thm-D.D} and \cite[Remark~III.4.10.9(b)]{Ama95a}.
\end{proof}
\begin{cor}\label{cor-D.A}
There exist
\hb{\om>0} such that
$$
\lda+\cA
\in\Lis\bigl(W_q^{2,\lda}(M;\rho),L_q^{\lda-2}(M;\rho)\bigr),
\ \lda\geq\om.
$$
\end{cor}
\begin{proof}
\cite[Theorem~I.2.2]{Ama95a}.
\end{proof}
\begin{rem}\label{rem-D.D}
Theorem~\ref{thm-D.D} is a maximal Sobolev space regularity statement.
Such results form the basis for proving existence results for
quasilinear parabolic equations on singular manifolds.
In this connection the minimal regularity assumptions~\eqref{D.aa}
are crucial.

\smallskip
Theorem~1.23 in~\cite{Ama17b} also contains a maximal H\"older space
regularity theorem on urR manifolds.
Y.~Shao and G.~Simonett~\cite{ShS14a} have implemented the Da~Prato--Grisvard
theorem~\cite{DaPG79a} on continuous maximal regularity on urR manifolds and
applied it to the Yamabe flow (also see \cite{LSS20a},~\cite{LS16a}).
Clearly, building on these results, similarly as we proved
Theorem~\ref{thm-D.D}, we can establish maximal H\"older space
and continuous regularity theorems in weighted spaces.
\qed
\end{rem}
It should be observed that the operator~\eqref{D.A} is non-degenerate. Linear
parabolic problems with degenerate coefficients have been investigated
in~\cite{Ama16a}, \cite{Ama20a}, \cite{Ama21a}, and nonlinear ones by Y.~Shao
\cite{Sha16b}, \cite{Sha16a},~\cite{Sha18a} and in~\cite{Ama20b}.

\smallskip
The fact that the weighted Sobolev spaces on~$(M,g,\rho)$ can be expressed
in terms of Sobolev spaces on~$\Mwhg$  has also been noted in
\cite{AGNi19a} and in M.~Kohr and V.~Nistor~\cite{KoN22a}. Based on this
observation, the authors
apply earlier results (e.g.,~\cite{AGN19c},~\cite{GroNi20a}) to derive
regularity and isomorphism theorems for elliptic boundary value problems.
\section{Model Cusps}\label{sec-M}
As mentioned in the introduction, it remains to exhibit concrete classes of
singular manifolds. This is done in this and the following two sections.

\smallskip
We begin with very simple but important one-dimensional urR manifolds.
Hereafter,
$$
\bt\quad
I:=(0,1]\text{ and $R$ is a cusp characteristic},
$$
that is,
\begin{equation}\label{M.IR}
\bal
{\rm(i)}\quad
&R\in C^\iy\bigl(I,(0,1]\bigr),\ {\textstyle\int_I}dt/R(t)=\iy.\cr
{\rm(ii)}\quad
&\|R^{j-1}\pl^jR\|_{BC(I)}\leq c(j),\ j\in\thBN.
\eal
\end{equation}
\begin{exs}\label{exs-M.C}
(a)
\uti{Power characteristics}
\hb{\sR_\al(t):=t^\al},
\ \hb{\al\geq1}.

\smallskip
(b)
\uti{Exponential characteristics}
\hb{\exp_{[\al,\ba]}(t):=e^{\al(1-t^{-\ba})}} for
\hb{\al,\ba>0}.\qed
\end{exs}
\begin{prop}\label{prop-M.I}  $(I,\ dr^2/R^2)$ is a urR manifold.
\end{prop}
\begin{proof}
We set
\hb{\rho(s):=\int_s^1dt/R(t)} for
\hb{s\in I}. Then
\hb{\rho\in\Diff(I,\BR_+)} and
$$
(\rho^*ds)(r)=d\rho(r)=-dr/R(r).
$$
Thus
\hb{\rho^*(\BR_+,ds^2)=(I,\ dr^2/R^2)}. This shows that $\rho$~is an
isometric diffeomorphism from
\hb{(I,\ dr^2/R^2)} onto~$(\BR_+,ds^2)$.
Examples~\ref{exs-U.S}(a) and~(d) yield the claim.
\end{proof}
Let
\hb{1\leq m\leq\gsm}. Assume that either
\begin{equation}\label{M.S}
\bal
&B\text{ is a compact \hb{(m-1)}-dimensional
 submanifold}\cr
\noalign{\vskip-1\jot}
&\text{of the unit sphere $\BS^{\gsm-1}$ in $\BR^\gsm$ and}\cr
&\qquad\qquad C(B):=\{\,tb\ ;\ t\in I,\ b\in B\,\}\is\BR^\gsm,
\eal
\end{equation}
or
\begin{equation}\label{M.R}
\bal
&B\text{ is a compact \hb{(m-1)}-dimensional
 submanifold of $\BR^{\gsm-1}$ and}\cr
&K(R,B):=\bigl\{\,\bigl(t,R(t)b\bigr)
 \ ;\ t\in I,\ b\in B\,\bigr\}\is\BR\times\BR^{\gsm-1}=\BR^\gsm.
\eal
\end{equation}
Here $tb$~is identified with~$t\ia_B(b)$, where $\ia_B$~denotes the inclusion
 map
 \hb{B\hr\BR^\gsm} if \eqref{M.S} applies, resp.\
 \hb{B\hr\BR^{\gsm-1}} otherwise. For the sake of a uniform presentation
 we write
$$
\bal
&Z=Z(R,B)\text{ for either $C(B)$ or }K(R,B),\cr
\noalign{\vskip-1\jot}
&\text{where it is understood that \hb{R=\sR_1} in the first case}.
\eal
$$
It is said to be~a \emph{smooth model cusp} (\emph{cone} if
\hb{Z=C}) \emph{in}~$\BR^\gsm$. `Smooth' expresses the fact that the
\emph{base}~$B$ is a smooth manifold.
\begin{rem}\label{rem-M.B}
It is not assumed that $B$~be connected. For examle, if
\hb{m=1} and
\hb{\gsm=3}, then $B$~is a \hbox{$0$-dimensional} submanifold of~$\BS^2$ if
\hb{Z=C}, resp.\ of~$\BR^2$ otherwise. Hence
\hb{B=\{b_1,\ldots,b_k\}} for some
\hb{k\in\thBN}. If
\hb{Z=C}, then $Z$~consists of $k$~straight lines of length~$1$ in~$\BR^3$
emanating from the origin, but not containing it. Assume that
\hb{Z=K}. Then $Z$~is a pair-wise disjoint union of $k$~smooth curves of
finite length in~$\BR^3$, which originate from
\hb{0\in\BR^3} also. Given any two
\hb{a,b\in B}, write~$E$ for the plane containing $0$, $a$,~$b$. It
encompasses the two curves with the endpoints $a$ and~$b$. If
\hb{\al=1}, then they approach the origin transversally, that is,
they form an angle in~$E$. Otherwise, they determine a cusp in~$E$,
that is, $Z$~represents~a `bouquet of flowers'.
\qed
\end{rem}
The map
$$
f_Z\colon I\times B\ra Z
\quad (t,b)\mt
\begin{cases}
tb,
&\text{if }Z=C,\cr
\bigl(t,R(t)b\bigr),
&\text{otherwise},
\end{cases}
$$
is a diffeomorphism, the \emph{stretching diffeomorphism}. Note that
\hb{I\times B} is a manifold with corners if
\hb{\pl B\neq\es}. Hence
\hb{Z=f_Z(I\times B)} is such a manifold also. It is easy to see that
everything established from Section~\ref{sec-U} onwards extends naturally
to such corner manifolds. (In this case
\hb{Q_\ka^m:=[0,1)^2\times(-1,1)^{m-2}} has to be added
to~\eqref{U.Q} if $U_\ka$~is a neighborhood of a corner point.)

\smallskip
On~$B$ we introduce the pull-back metric
\hb{g_B:=\ia_B^*g_\ell}, where
\hb{\ell=\gsm} if
\hb{Z=C}, and
\hb{\ell=\gsm-1} otherwise. We endow~$Z$ with the metric
\begin{equation}\label{M.gZ}
g_Z:=f_{Z\,*}(dt^2+R^2g_B).
\end{equation}
It is equivalent to the metric induced by the embedding
\hb{\ia_Z\colon Z\hr\BR^\gsm}:
\begin{lem}\label{lem-M.fg}
\hb{g_Z\sim\ia_Z^*g_\gsm}.
\end{lem}
\begin{proof}
We identify~$f_Z$ with
\hb{\ia_Z\circ f_Z}.

\smallskip
(1)
Suppose that
\hb{Z=K}. Let $\gK$ be a ur atlas for~$B$ and put
$$
h_\ka:=\ia_B\circ\ka^{-1}\in C^\iy(Q_\ka^{m-1},\BR^{\gsm-1}),
\ \ka\in\gK.
$$
Then
$$
f_\ka:=\bigl((t,y)\mt (t,R(t)h_\ka(y))\bigr)
\in C^\iy(I\times Q_\ka^{m-1},\,\BR^\gsm)
$$
is a local parametrization of~$Z$. From now on,
the indices $i$ and~$\al$ run from $2$ to~$\gsm$ and from $2$ to~$m$,
respectively. Then
$$
df_\ka^1=dt,
\quad df_\ka^i=\dot Rh_\ka^idt+R\,\pl_\al h_\ka^idy^\al.
$$
Hence, writing
\hb{\vsdot=\vsdot_\gsm},
$$
\bal
f_\ka^*g_\gsm
&=dt^2+\sum_i|\dot R h_\ka^idt+R\,\pl_\al h_\ka^idy^\al|^2\cr
&=(1+\dot R^2\,|h_\ka|^2)dt^2+2R\dot R dt
 \sum_ih_\ka^i\,\pl_\al h_\ka^idy^\al+R^2\sum_i|\pl_\al h_\ka^idy^\al|^2.
\eal
$$
Since $B$~is compact, we can fix
\hb{\da\geq1} such that
\hb{|h_\ka|\leq\da} for
\hb{\ka\in\gK}. The Cauchy--Schwarz inequality yields
$$
\Big|2R\dot R dt\sum_ih_\ka^i\,\pl_\al h^idy^\al\Big|
\leq\ve^{-1}\dot R^2\,|h_\ka|^2\,dt^2
+\ve R^2\sum_i|\pl_\al h_\ka^idy^\al|^2\ .
$$
This implies that $f_\ka^*g_\gsm$~can be estimated from below by
$$
\bigl(1+(1-\ve^{-1})\dot R^2\,|h_\ka|^2\bigr)dt^2
+(1-\ve)R^2\sum_i|\pl_\al h_\ka^idy^\al|^2.
$$
We fix
\hb{\ve\in(0,1)} close to~$1$ such that
\hb{(\ve^{-1}-1)\,\|\dot R\|_\iy^2\,\da^2<1}. This guarantees, due to
\hb{\ka_*g_B=h_\ka^*g_{\gsm-1}}, that
$$
f_\ka^*g_\gsm\geq \Bigl(dt^2+R^2\sum_i|\pl_\al h_\ka^idy^\al|^2\Bigr)\Big/c
=(dt^2+R^2\ka_*g_B)/c.
$$
A~similar argument with
\hb{\ve=1} yields
$$
f_\ka^*g_\gsm\leq c(dt^2+R^2\ka_*g)
$$
\hbox{$\ka$-uniformly}. The claim follows.

\smallskip
(2)
Assume
\hb{Z=C}. Then
\hb{f_\ka(t,y)=th_\ka(y)} for
\hb{y\in Q_\ka^{m-1}}, where now
\hb{h_\ka(y)\in\BS^{\gsm-1}\hr\BR^\gsm}. Hence
\hb{|h_\ka|=1}. Consequently, letting $j$ run from~$1$ to~$\gsm$
and $\ba$~from~$1$ to~$m$,
$$
2\sum_jh_\ka^j\,\pl_\ba h_\ka^j=\sum_j\pl_\ba(h_\ka^j)^2=0.
$$
From this,
\hb{\pl_tf_\ka=h_\ka}, and
\hb{\pl_\ba f_\ka=t\,\pl_\ba h_\ka} we get
\begin{equation}\label{M.fg}
f_\ka^*g_\gsm=dt^2+t^2\ka_*g_B,
\ \ka\in\gK.
\end{equation}
Thus the assertion applies in this case also.
\end{proof}
Note that, by \eqref{M.fg},
\begin{equation}\label{M.gC}
g_C=\ia_C^*g_\gsm.
\end{equation}
We define
\hb{r_Z\in C^\iy\bigl(Z,(0,1]\bigr)} by
\begin{equation}\label{M.rZ}
r_Z(x):=
\begin{cases} |x|_\gsm,
&\text{if }Z=C,\cr
R(x^1),
&\text{otherwise}.\cr
\end{cases}
\end{equation}
\begin{thm}\label{thm-M.g}
Set
\hb{\wh{g}_Z:=g_Z/r_Z^2}. Then $(M,\wh{g}_Z)$ is a urR manifold.
\end{thm}
\begin{proof}
Since
\hb{r_Z=f_{Z\,*}R=R\circ f_Z^{-1}},
$$
\wh{g}_Z=(f_{Z\,*}R)^{-2}f_{Z\,*}(dt^2+R^2g_B)=f_{Z\,*}(R^{-2}dt^2+g_B).
$$
Thus the claim follows from Proposition~\ref{prop-M.I} and
Examples \ref{exs-U.S}\hbox{(b)--(d)}.
\end{proof}
\begin{thm}\label{thm-M.r}
$r_Z$~is a singularity function for~$(Z,g_Z)$.
\end{thm}
\begin{proof}
(1)
We claim that
\begin{equation}\label{M.Rk}
\|R^k\pl^k(\log R)\|_{BC(I)}\leq c(\al,k),
\ k\in\thBN.
\end{equation}
We proceed by induction.

\smallskip
Set
\hb{a:=\pl(\log R)} so that
\hb{R a=\pl R}. Hence
\hb{|R\pl(\log R)|\leq c(0)}. Thus the assertion is proved if
\hb{k=0}.

\smallskip
Assume
\hb{k\in\thBN} and
\begin{equation}\label{M.i}
|R^{i+1}\pl^ia|\leq c(\al,i),
\ 0\leq i\leq k-1.
\end{equation}
Then
$$
\bal
R^{k+1}\pl^ka
&=R^k\pl^k(R a)
 -\sum_{i=0}^{k-1}\binom kiR^k\pl^i R\,\pl^{k-i}a\cr
&=R^k\pl^{k+1}R
 -\sum_{i=1}^{k-1}\binom ki(R^{i-1}\pl^iR)(R^{k-i+1}\pl^{k-i}a)\ .
\eal
$$
Using this and \eqref{M.IR}(ii), we deduce that \eqref{M.i} holds also for
\hb{i=k}. This proves the claim.

\smallskip
(2)
Set
\hb{g_0:=dt^2+R^2g_B}. Then
$$
\na_{g_0}\om=\pl_t\om_1\oplus\na_{g_B}\om_2,
\ \om=\om_1\oplus\om_2\in T^*I\oplus T^*B=T^*(I\times B).
$$
Hence
\hb{\na_{g_0}^k\om=\pl_t^k\om_1} if
\hb{\om_2=0}. Consequently,
$$
\bal
&\big\|r_Z^{k+1}\,|\na_{g_Z}^kd(\log r_Z)|_{(g_Z)_0^{k+1}}
 \,\big\|_{BC(Z)}\cr
&\qquad\qquad{}
=\big\|f_{Z\,*}\bigl(R^{k+1}\,|\pl_t^{k+1}(\log R)|\bigr)
 \,\big\|_{BC(Z)}\cr
&\qquad\qquad{}
 =\|R^{k+1}\pl_t^{k+1}(\log R)\|_{BC(I)}\leq c(\al,k)
\eal
$$
for
\hb{k\in\BN} by step~(1). Due to Theorem~\ref{thm-M.g}, the assertion follows.
\end{proof}
\section{Manifolds with Point Singularities}\label{sec-P}
We denote, for
\hb{m\in\thBN}, by~$\BB^m$ the open unit ball in
\hb{\BR^m=\Rmgm}.

\smallskip
Let $S$ be a \hbox{$0$-dimensional} submanifold of~$\BR^\gsm$. This means
that each
\hb{p\in S} has an open neighborhood~$\cU_p$ such that
\hb{\cU_p\cap S=\{p\}}. The uniform regularity of~$(\BR^\gsm,g_\gsm)$ and the
separability of~$\BR^\gsm$ imply that
$$
\inf\{\,|p-q|_\gsm\ ;\ p,q\in\cS,\ p\neq q\,\}>0.
$$
This shows that
part~(iv) of the subsequent definition is meaningful.
\begin{defn}\label{def-P.D}
Let
\hb{1\leq m\leq\gsm}. Assume:
\begin{itemize}
\item[\rm(i)]
$\Mg$~is an \hbox{$m$-dimensional} Riemannian submanifold
of $(\BR^\gsm,g_\gsm)$.
\item[\rm(ii)]
\hb{\cS=\cS(M):=\cl_{\BR^\gsm}(M)\setminus M} is a
\hbox{$0$-dimensional} submanifold of~$\BR^\gsm$, the
\emph{singularity set of}~$M$.
\item[\rm(iii)]
$\cZ$~is a finite set of \hbox{$m$-dimensional} model cusps
\hb{Z=Z(R,B)} in~$\BR^\gsm$, and
\hb{Z(r):=Z\cap r\BB^\gsm} for
\hb{0<r<1}.
\item[\rm(iv)] There exist
\hb{\ve\in(0,1)}, for each
\hb{p\in\cS} a~model cusp
\hb{Z_p\in\cZ}, an open neighborhood~$\wt{U}_p$ of~$\ol{Z_p(\ve)}$ in~$\BR^\gsm$,
and an injective immersion~$\wt{\psi}_p$ from~$\wt{U}_p$
into~$\BR^\gsm$ satisfying
\hb{\wt{\psi}_p(0)=p},
$$
M(p,\ve):=\wt{\psi}_p\bigl(Z_p(\ve)\bigr)\is M,
$$
\hb{M(p,\ve)\cap M(q,\ve)=\es} for
\hb{p,q\in M} with
\hb{p\neq q}, and,
setting
\hb{\psi_p:=\wt{\psi}_p|Z_p},
\begin{equation}\label{P.p}
\|\psi_p\|_{k,\iy}\leq c(k),
\ p\in\cS,
\ k\in\thBN.
\end{equation}
Then $M(p,\ve)$ is called \emph{cusp of~$M$ at}~$p$ (although
\hb{p\notin M}) \emph{of type}~$Z_p$, and
$$
M(\cS,\ve):=\bigcup_{p\in\cS}M(p,\ve)
$$
is a \emph{cuspidal neighborhood of~$M$ along~$\cS$ of type~$\cZ$}. Also,
$\Mg$~is a \emph{(singular) manifold with (smooth) point
singularities of type}~$\cZ$ if it possesses a cuspidal neighborhood
along~$\cS$ of type~$\cZ$.
\end{itemize}
Of course, $M(p,\ve)$ is called \emph{cone} at~$p$
if its model cusp is a cone.\qed
\end{defn}
We introduce a \emph{cuspidal chart}~$\vp_p$ for~$M$ at~$p$ by
$$
\vp_p:=\psi_p^{-1}|M(p,\ve)\colon M(p,\ve)\ra Z_p(\ve),
$$
which is a diffeomorphism, and
\hb{\{\,\vp_p\ ;\ p\in\cS\,\}} is a \emph{cuspidal atlas} for
$M(\cS,\ve)$. Henceforth,
\begin{equation}\label{P.Z}
\bal
&M\text{ is a manifold with point singularities and}\cr
\noalign{\vskip-1\jot}
&M(\cS,\ve)\text{ is a cuspidal neighborhood of type }\cZ.
\eal
\end{equation}
\begin{rems}\label{rems-P.M}
(a)
Let
\hb{p\in\cS} and
\hb{Z_p=Z(R,B)\in\cZ}. Then
\hb{\pl Z_p(\ve)\neq\es} iff
\hb{\pl B\neq\es}, and
\hb{\pl Z_p(\ve)=Z_p(R,\pl B)\cap\ve\BB^\gsm}. Thus
$$ \pl M(p,\ve)=\vp_p^{-1}\bigl(\pl Z_p(\ve)\bigr).
$$
This shows that $p$~is a singular point `at the boundary of~$M$'.

\smallskip
(b)
Assume that $p$~is a conical point of type~$C(B)$. Then
\hb{\pl B=\es} iff
\hb{B=\BS^{m-1}\hr\BS^{\gsm-1}}. Since
\hb{C(\BS^{m-1})=\BB^m\setminus\{0\}},
this means that `$M$~has a hole at~$p$'.

\smallskip
(c)
If $p$~is a conical point of type
\hb{C(\BS^{m-1}\cap\BH^m)}, then `$p$~is a hole in the boundary~$\pl M$'.

\smallskip
(d)
A~manifold with point singularities can have countably infinitely many cusps,
but of finitely many types only. In this case $M$~is unbounded.\qed
\end{rems}
Let $(N,h)$ be a Riemannian manifold and
\hb{A\is N}. An atlas~$\gK$ for~$N$ is \emph{ur on}~$A$ if
Definition~\ref{def-U.K} applies for
\hb{\ka\in\gK_A:=\{\,\nu\in\gK\ ;\ U_\nu\cap A\neq\es\}}. Moreover,
$(N,h)$~is \emph{ur on}~$A$ if Definition~\ref{def-U.gg} holds
with $g$~replaced by~$h$ and $\gK$ with~$\gK_A$. Note that
$\gK_A$~is a finite set if $A$~is compact and $\gK$~is ur on~$A$.

\smallskip
We say that $\Mg$ is ur \emph{off}~$\cS(M)$ if $\Mg$ is ur on
$$
M(\cS,r)^c:=M\setminus M(\cS,r)
$$
for each
\hb{r\in(0,\ve)}.
\begin{exs}\label{exs-P.S}
(a)
A~nonempty boundary~$\pl M$ is said to be \emph{almost regularly embedded} if
\hb{m=\gsm} and
\hb{\pl M\cap M(\cS,r)^c} has for each
\hb{r\in(0,\ve)} a~uniform geodesic collar. If this prevails, then $\Mg$
is ur off~$\cS(M)$.
\begin{proof}
This follows from the fact that there are ur atlases for~$\BR^m$
whose coordinate patches have arbitrarily small diameters.
\end{proof}

\smallskip
(b)
Let $M$ be bounded in~$\BR^\gsm$. Then $\Mg$~is ur off~$\cS(M)$.
\begin{proof}
Since $M(\cS,r)^c$ is compact, the localized version of
Example~\ref{exs-U.S}(b) yields the assertion.
\end{proof}

\smallskip
(c)
If
\hb{m=\gsm} and $\pl M$~is relatively compact, then $\Mg$ is ur
off~$\cS(M)$.
\begin{proof}
In this situation $\pl M$~is almost regularly embedded.
\end{proof}
\end{exs}
We fix
\hb{0<\ve_0<\ve_1<\ve} and an increasing function
\hb{\om\in C^\iy[0,\ve]} satisfying
\hb{\om(t)=0} if
\hb{t\leq\ve_0}, and
\hb{\om(t)=1} for
\hb{t\geq\ve_1}. Then
$$
\om_p:=\vp_p^*\om\in C^\iy\bigl(M(p,\ve),(0,1]\bigr).
$$
Hence
\begin{equation}\label{P.rp}
\rho_p:=(1-\om_p)\vp_p^*r_{Z_p}+\om_p
\in C^\iy\bigl(M(p,\ve),(0,1]\bigr),
\ p\in\cS.
\end{equation}
The \emph{cusp characteristic}~$\rho$ for~$M$, defined by
\begin{equation}\label{P.r}
\rho:=
\begin{cases} \rho_p\text{ on }M(p,\ve),
&p\in\cS,\cr
1\text{ on }M(\cS,\ve)^c,
\end{cases}
\end{equation}
belongs to $C^\iy\bigl(M,(0,1]\bigr)$. We also introduce a Riemannian
metric~$\ol{g}$ on~$M$ by
\begin{equation}\label{P.gp}
\ol{g}_p:=(1-\om_p)\vp_p^*g_{Z_p}+\om g,
\ p\in\cS,
\end{equation}
and
\begin{equation}\label{P.g}
\ol{g}:=
\begin{cases}
\ol{g}_p\text{ on }M(p,\ve),
&p\in\cS,\cr
g\text{ on }M(\cS,\ve)^c.
\end{cases}
\end{equation}
\begin{thm}\label{thm-P.M}
If $\Mg$ is ur off~$\cS(M)$ and
\hb{\wh{g}:=\ol{g}/\rho^2}, then $\Mwhg$ is ur.
\end{thm}
\begin{proof}
Set
\hb{M(p,r):=\vp_p^{-1}\bigl(Z_p(r)\bigr)},
\ \hb{0<r<\ve}, and
\hb{M(p,r]:=\cl_M\bigl(M(p,r)\bigr)}.

\smallskip
(1)
Theorem~\ref{thm-M.g} guarantees that $(Z_p,\wh{g}_{Z_p})$ is ur.
Thus it is ur on
$$
Z_p(\ve_1]:=Z_p\cap\ve_1\ol{\BB}\vph{\BB}^\gsm.
$$
Example~\ref{exs-U.S}(d) yields that
$\bigl(M(p,\ve),\vp_p^*\wh{g}_{Z_p}\bigr)$ is ur on~$M(p,\ve_1]$. By the
finiteness of $\cZ$ and \eqref{P.p} this holds uniformly w.r.t.\
\hb{p\in\cS}.

\smallskip
(2)
Note that
\begin{equation}\label{P.rg}
\rho_p=\vp_p^*r_{Z_p},
\ \ol{g}_p=\vp_p^*g_{Z_p}
\quad\text{on }Z_p(\ve_0].
\end{equation}
Thus
\hb{\vp_p^*\wh{g}_{Z_p}=\wh{g}} on $M(p,\ve_0]$.

\smallskip
(3)
Since $\Mg$~is ur off~$\cS(M)$, \,$\Mg$~is ur on $M(\cS,\ve_0)^c$. Note that
\hb{A_p:=M(p,\ve_1]\cap M(p,\ve_0)^c} is compact. Let $\ol{\gK}$ be an atlas
for~$M$ which is ur on~$M(\cS,\ve_0)^c$. Its restriction to $M(p,\ve)$
is an atlas which is ur on~$A_p$. Hence
\hb{\ol{\gK}_p:=\ol{\gK}_{A_p}} is finite. Similarly, choose an
atlas~$\gZ_p$ for~$Z_p(\ve)$ which is ur on~$Z_p(\ve_1]$. Then
$\vp_p^*\gZ_p$~is an atlas for~$M(p,\ve)$ which is ur on~$M(p,\ve_1]$.
Consequently,
\hb{\wh{\gK}_p:=(\vp_p^*\gZ_p)_{A_p}} is finite as well. Thus
\hb{\ol{\gK}_p\cup\wh{\gK}_p} is an atlas for the neighborhood
$$
\bigcup_{\ol{\ka}\in\ol{\gK}_p}U_{\ol{\ka}}
\cup\bigcup_{\wh{\ka}\in\wh{\gK}_p}U_{\wh{\ka}}
$$
of~$A_p$ which is finite and ur on~$A_p$. From this we infer that
$\bigl(M(p,\ve),\ol{g}_p\bigr)$ is ur on~$A_p$. Since
\hb{\min\bigl\{\,\rho_p(q)\ ;\ q\in A_p\,\bigr\}>0}, it follows that
\hb{\wh{g}\sim\ol{g}_p} on~$A_p$. Consequently, $\Mwhg$
is ur on~$A_p$. The finiteness of $\cZ$ and \eqref{P.p} guarantee that this
holds uniformly w.r.t.\
\hb{p\in\cS}. From this and
\hbox{\eqref{P.rp}--\eqref{P.g}} the assertion follows.
\end{proof}
\begin{rem}\label{rem-P.M}
It is seen that a different choice of~$\ve$, the cut-off function~$\om$, and
of an equivalent cuspidal atlas leads to a ur metric for~$M$ which is
equivalent to~$\wh{g}$. Thus this theorem means that we `uniformly
regularize' the singular manifold~$\Mg$ by means of a ur metric~$\wh{g}$
which differs from~$g$ only arbitrarily close to the singularity
set~$\cS(M)$.
\qed
\end{rem}
Let $(M_i,g_i)$ be Riemannian manifolds and
\hb{f\colon(M_1,g_1)\ra(M_2,g_2)} be an isometric diffeomorphism. Set
\hb{\na_i:=\na_{g_i}} and define
$$
f^*\na_2\text{ by }
(f^*\na_2)\om=f^*\bigl(\na_2(f_*\om)\bigr),
\ \om\in C^\iy(T^*M_1).
$$
Then
\hb{f^*\na_2=\na_1} (see \cite[Theorem~X.2.3.2]{Ama25a}). Hence
\begin{equation}\label{P.fN}
f^*\na_2^k=(f^*\na_2)^k=\na_1^k.
\end{equation}
\begin{thm}\label{thm-P.f}
$\rho$~is a singularity function for~$(M,\ol{g})$.
\end{thm}
\begin{proof}
Due to Theorem~\ref{thm-P.M}, it remains to prove that
\begin{equation}\label{P.log}
\|\rho^{k+1}\na_{\ol{g}}^kd(\log\rho)\|_{BC(M)}\leq c(k),
\ k\in\BN.
\end{equation}
On $M(\cS,\ve_0)^c$,
\hb{\ul{\rho}\leq\rho\leq1} for some
\hb{\ul{\rho}>0}.  Therefore \eqref{P.log} is trivially true on
$M(\cS,\ve_0)^c$. Hence
it remains to prove that \eqref{P.log} holds on $M(\cS,\ve_0)$.
Note that, by \eqref{P.rg} and \eqref{P.fN},
$$
\vp_{p\,*}\bigl(\rho^{k+1}\na_{\ol{g}_p}^kd(\log\rho)\bigr)
=r_{Z_p}^{k+1}\na_{g_{Z_p}}^kd(\log r_{Z_p})
\text{ on }Z_p(\ve_0].
$$
Thus the assertion follows from Theorem~\ref{thm-M.r} and the
finiteness of~$\cZ$.
\end{proof}
Manifolds with point singularities are the simplest class of singular
manifolds. Larger families of uniformly regularizable singular manifolds
are presented in \cite{Ama15a} and in~\cite{Ama25a}. They comprise manifolds
with cuspidal corners, cuspidal wedges, intruding cones, etc.

\smallskip
The starting point for a thorough analysis of singular manifolds and of
differential operators thereon is the paper by
V.A.~Kondratiev~\cite{Kon67a}. Since then, there has appeared an
inextricable flood of publications and it is impossible to do justice to the
authors. We simply mention the prolific works of V.G. Maz{\cprime}ya and
\hbox{B.-W.} Schulze. Together with coworkers, the first author developed
further Kondratiev's approach (see~\cite{KMR97a} for an early influential
presentation). The main interest of the second author concerns algebras of
pseudo-differential operators on manifolds with singularities. (Also the work by R.B.~Melrose~\cite{Mel93a} should be mentioned). Schulze builds on the
Mellin transform and cone differential operators of Fuchs-type (see
\cite{Schu91a}, or the more recent
book by V.E. Nazaikinskii, A.Yu. Savin, \hbox{B.-W.} Schulze, and
B.Yu. Sternin~\cite{NSSS06a} for accessible accounts). Both directions of
research focus on singular function expansions which are not possible by our
technique. E.~Schrohe and coauthors develop a maximal
regularity analysis of parabolic evolution equations on manifolds with
conical singularities. They employ Schulze's Mellin--Sobolev spaces and
implement the functional analytic Dore--Venni theorem (e.g.,~\cite{Ama95a}).
This
requires the semigroup generator to have bounded imaginary powers, which is
established for a class of cone differential operators containing the
Laplace--Beltrami operator; see N.~Roidos and E.~Schrohe~\cite{RoS14a},
E.~Schrohe and J.~Seiler \cite{SchS18a},
S.~Coriasco, E.~Schrohe, and J.~Seiler \cite{CSS03a},~\cite{CSS07a}.
In contrast, Theorem~1.23 in~\cite{Ama17b}, which is the fundament for
Theorem~\ref{thm-D.D}, is based on a Fourier multiplier theorem on~$\BR^m$
and on the $(\cR,\cR^c)$ localization technique, which apply with the same
ease to parabolic equations of arbitrary order with values in spaces of
sections of vector bundles.
\section{Conical Singularities}\label{sec-K}
Assume that
\hb{Z_p=C(B_p)} with
\hb{p\in\cS}. If
\hb{q\in M(p,\ve)}, we write
\hb{d_p(q)=d_p^M(q)} for the infimum of the lengths of all smooth curves
\hb{\ga\colon[0,1]\ra\BR^\gsm} satisfying
\hb{\ga(0)=p},
\ \hb{\ga(1)=q}, and
\hb{\ga(0,1]\is M(p,\ve)}. Thus $d_p^M(q)$~is the Riemannian distance
in~$M(p,\ve)$ from~$q$ to the `point~$p$ at infinity' of~$\Mwhg$.
\begin{thm}\label{thm-K.M}
Let $\Mg$ be an
\hbox{$m$-dimensional} Riemannian submanifold of $(\BR^\gsm,g_\gsm)$ with
conical singularities. Select a {\rm distance function
\hb{\da\in C^\iy(M)} for}~$\cS(M)$, that is,
$$
\da\sim d_p\text{ on }M(p,\ve),
\ p\in\cS,
\quad \da=1\text{ on }M(\cS,\ve)^c.
$$
Then
$$
u\mt\sum_{j=0}^k\big\|\da^{-\lda+j-m/q}\ |\na^ju|_{g_0^j}\,\big\|_{L_q(M)},
\ 1\leq q\leq\iy,
$$
is a norm for $W_q^{k,\lda}(M;\rho)$ if
\hb{q<\iy}, resp.\ for $BC^{k,\lda}(M;\rho)$ if
\hb{q=\iy}.
\end{thm}
\begin{proof}
Let $C_p$ be the model cone for~$M(p,\ve)$ and assume that $\vp_p$~is a
conical chart. We infer from \eqref{M.gC} and~\eqref{M.rZ} that
\hb{r_{C_p}=\vsdot_\gsm}. Hence
\hb{d_0^{C_p}=r_{C_p}}. Consequently,
\hb{\vp_p^*r_{C_p}=\vp_p^*d_0^{C_p}\sim d_p^M}. Now the claim
follows from~\eqref{P.rp}.
\end{proof}
In order to not overstretch the present paper, we consider smooth
singularities only. However, manifolds with corners, edges, interfaces, cuts,
etc.\ play an important role, in numerical analysis in particular. In such
situations it is also possible---though much more technical---to carry out a
uniform regularization along the lines of the proof of Theorem~\ref{thm-P.M}
(see~\cite{Ama25a}).

\smallskip
A~different approach is due to
C.~B\u{a}cu\c{t}\u{a}, A.L. Mazzucato, V.~Nistor, and
L.~Zikatanov~\cite{BMNZ10a}. These authors employ~a `desingularization'
technique for curvilinear polyhedra based on the theory of Lie manifolds
introduced by B.~Ammann, R.~Lauter, and V.~Nistor~\cite{ALN04a} and
B.~Ammann, A.D. Ionescu, and V.~Nistor~\cite{AIN06a}.
It is essentially equivalent to our uniform regularization technique. This
implies that the general theory of weighted spaces, exposed in Sections
\ref{sec-W} and~\ref{sec-D}, applies to general bounded curvilinear polyhedra
and singularity functions which are equivalent to the distance to the full
singularity set~$\cS(M)$.

\smallskip
Lastly, we consider the important classical Euclidean setting. Precisely,
we suppose that
\begin{equation}\label{K.m}
\bal
&\Mg\text{ is an \hbox{$m$-dimensional} Riemannian submanifold}\cr
\noalign{\vskip-1\jot}
&\text{of $(\BR^m,g_m)$ with conical singularities}\cr
\noalign{\vskip-1\jot}
&\text{such that }\pl M\text{ is almost regularly embedded}.
\eal
\end{equation}
The last condition is satisfied (see Example~\ref{exs-P.S}(c)) if $\pl M$~is
relatively compact, thus, in particular, if $M$~is bounded, which is the case
most often looked at in the literature.

\smallskip
Let $\da$ be a distance function for~$\cS(M)$,
\ \hb{k\in\BN},
\ \hb{1\leq q\leq\iy}, and
\hb{a\in\BR}. The \emph{Kondratiev space}~$K_{q,a}^k(M)$ consists of all
\hb{u\in L_{1,\loc}(M)} whose distributional derivatives (in~$\ci M$)
of order at most~$k$ are such that
$$
\|u\|_{K_{q,a}^k}
:=\sum_{|\al|\leq k}\|\da^{|\al|-a}\pl^\al u\|_{L_q(M)}<\iy.
$$
It is endowed with this norm.
\begin{thm}\label{thm-K.K}
Let \eqref{K.m} be satisfied. Then
$$
K_{q,a}^k(M)\doteq
\begin{cases}
W_q^{k,\,a-m/q}(M;\rho),
&\text{if }q<\iy,\cr
BC^{k,a}(M;\rho),
&\text{otherwise}.
\end{cases}
$$
\end{thm}
\begin{proof}
Now
\hb{g=g_m}. Thus
\hb{\na^ju=\pl^ju\in C(M,\cL^j)} for
\hb{u\in C^j(M)}, where $\cL^j$~is the space of \hbox{$j$-linear} functions
on~$\BR^m$. It is given the usual norm~%
\hb{\vsdot_{\cL_j}}. Set
\hb{\lda:=a-m/q}. Then
$$
\bal
\sum_{j=0}^k\big\|\da^{-\lda+j-m/q}\ |\na^ju|_{g_0^j}\,\big\|_{L_q(M)}
&\sim\sum_{j=0}^k\big\|\da^{j-a}\ |\pl^ju|_{\cL^j}\,\big\|_{L_q(M)}\cr
&\sim\sum_{|\al|\leq k}\|\da^{|\al|-k}\pl^\al u\|_{L_q(M)}.
\eal
$$
Since the distributional, that is weak, derivatives (on~$\ci M$) coincide
with the strong ones, and $\pl M$~is an
\hbox{$m$-dimensional} Lebesgue \hbox{$0$-set}, well-known arguments
yield that $K_{q,a}^k(M)$~is the completion of
\hb{(\cD(M),\Vsdot_{K_{q,a}^k})} in~$L_{1,\loc}(M)$ if
\hb{q<\iy}, resp.\ of
\hb{(BC^\iy(M),\Vsdot_{K_{\iy,a}^k})} in~$C(M)$ otherwise.
Now the assertion is clear.
\end{proof}
Kondratiev spaces on subdomains of Euclidean spaces have recently gained some
attraction, see St. Dahlke, M.~Hansen, C.~Schneider, and W.~Si\-ckel
\cite{DaHSS19a},~\cite{DahHSS20a},
St. Dahlke and C.~Schneider \cite{DaSch19a},~\cite{DaSch21a}, for example.
(The references in these papers to earlier works of various other authors
should also be taken into account.) These publications contain proofs of
embedding, interpolation, and point-wise multiplier theorems for a subclass
of the curvilinear polyhedra considered by V.~Nistor and coauthors. Some
investigations of Kondratiev spaces of fractional order are found
in~\cite{HaSch24a}.

\smallskip
Due to Theorem~\ref{thm-K.K} and the remarks following Theorem~\ref{thm-K.M},
the full theory exposed in Section~\ref{sec-W} applies to the clasical
Kondratiev spaces on curvilinear polyhedra.
\section*{Appendix: Notations and Conventions}\label{sec-A}
\renewcommand{\theequation}{\rm A.\arabic{equation}}
\setcounter{equation}{0}
We employ standard notations some of which we recall below. As for manifolds,
we refer to~\cite{Die69b}. A~detailed exposition is given in~\cite{Ama25a}.

\smallskip
We use~$c$ to denote a generic constant
\hb{\geq1}, whose value may be different in different formulas
but is always independent of the free variables in a given setting.
Real vector spaces are considered throughout.
The complex case can be covered by complexification.

\smallskip
Let $S$ be a nonempty set. On the vector space
\hb{\Ga(S):=\BR^S} of all real-valued functions on~$S$, an equivalence
relation~%
\hbox{${}\sim{}$} is defined by
\hb{f_1\sim f_2} iff
\hb{f_1/c\leq f_2\leq cf_1}. If $S$~is a subset of some vector space, then
\hb{\thS:=S\setminus\{0\}}.

\smallskip
Given Banach spaces $E$, $E_1$, and~$E_2$,
\ $\cL(E_1,E_2)$~is the Banach space of bounded linear maps from~$E_1$
into~$E_2$, and $\Lis(E_1,E_2)$ is the open subset of isomorphisms. The dual
of~$E$ is written~$E'$,
\ \hbox{$\pw=\pw_E\colon E'\times E\ra\BR$} is the duality pairing, and
\hb{a'\in\cL(E_2',E_1')} is the dual of
\hb{a\in\cL(E_1,E_2)}. As usual,
\hb{E_1\hr E_2} denotes continuous injection, and
\hbox{${}\sdh{}$} says that $E_1$ is also dense in~$E_2$. Furthermore,
\hb{E_1\doteq E_2} iff
\hb{E_1\hr E_2} and
\hb{E_2\hr E_1}, that is, $E_1$ and~$E_2$ are equal except for
equivalent norms.

\smallskip
Let $M$ be an \hbox{$m$-dimensional} manifold with (a~possibly empty)
boundary. We work in the smooth category and assume that the underlying
topological space is separable and metrizable. If $V$~is a vector bundle
over~$M$, then
\hb{\Ga(V)=\Ga(M,V)} is the $\Ga(M)$-module of sections of~$V$
(\emph{no} smoothness). Thus
\hb{\Ga(M)=\Ga(M\times\BR)}, where
\hb{M\times\BR} is the trivial bundle. If $V$~is a metric vector bundle, then
$C^\iy(V)$~is the $C^\iy(M)$~submodule of smooth sections.

\smallskip
As customary, $TM$ and~$T^*M$ are the tangent and cotangent bundles of~$M$.
Then
\hb{T_\tau^\sa M:=TM^{\otimes\sa}\otimes T^*M^{\otimes\tau}} is, for
\hb{\sa,\tau\in\BN}, the $(\sa,\tau)$-tensor bundle of~$M$, that is, the
vector bundle of all tensors on~$M$ being contravariant of order~$\sa$
and covariant of order~$\tau$. In particular,
\hb{T_0^1M=TM} and
\hb{T_1^0M=T^*M}, as well as
\hb{T_0^0M=M\times\BR}.

\smallskip
For
\hb{\nu\in\BN^\times} we put
\hb{\BJ_\nu:=\{1,\ldots,m\}^\nu}. Then, given local coordinates
\hb{\ka=(x^1,\ldots,x^m)} on an open subset~$U$ of~$M$ and setting
$$
\frac\pl{\pl x^{(i)}}
:=\frac\pl{\pl x^{i_1}}\otimes\cdots\otimes\frac\pl{\pl x^{i_\sa}},
\quad dx^{(j)}:=dx^{j_1}\otimes\cdots\otimes dx^{j_\tau}
$$
for
\hb{(i)=(i_1,\ldots,i_\sa)\in\BJ_\sa},
\ \hb{(j)\in\BJ_\tau}, the local representation of~a $(\sa,\tau)$-tensor
field
\hb{a\in\Ga(T_\tau^\sa M)} with respect to these coordinates is given by
$$
a=a_{(j)}^{(i)}\frac\pl{\pl x^{(i)}}\otimes dx^{(j)}
$$
with
\hb{a_{(j)}^{(i)}\in\Ga(U)}. We use the summation convention
for \hbox{(multi-\nobreak)}\linebreak[0]indices
labeling coordinates or
bases. Thus such a repeated index, which appears once as a superscript and
once as a subscript, implies summation over its whole range.

\smallskip
Suppose
\hb{1\leq s\leq\sa} and
\hb{1\leq t\leq\tau}. We write
\hb{(i)\in\BJ_\sa} in the form~$(i_s;i_{\wh{s}})$ where we set
\hb{(i_{\wh{s}}):=(i_1,\ldots,\wh{i_s},\ldots,i_\sa)\in\BJ_{\sa-1}}.
Then we define the
\emph{contraction map,~$\sC_t^s$, with respect to positions $s$ and}~$t$ by
$$
\sC_t^sa
:=a_{(k;j_{\wh{t}})}^{(k;i_{\wh{s}})}
\frac\pl{\pl x^{(i_{\wh{s}})}}\otimes dx^{(j_{\wh{t}})},
\ a\in\Ga(T_\tau^\sa M),
$$
where $k$~runs from~$1$ to~$m$. Hence
\hb{\sC_t^s\colon C^\iy(T_\tau^\sa M)
   \ra C^\iy(T_{\tau-1}^{\sa-1}M)} and this map is linear.

\smallskip
Let
\hb{0\leq\rho\leq\sa} and
$$
(a,b)\in\Ga(T_\tau^\sa M\oplus T_\rho^0 M)\ .
$$
We write
\hb{(i)\in\BJ_\sa} in the form~$(\ell)(k)$ with
\hb{(k)\in\BJ_\rho}. Then
$$
a=a_{(j)}^{(\ell)(k)}\frac\pl{\pl x^{(\ell)}}
\otimes\frac\pl{\pl x^{(k)}}\otimes dx^{(j)},
\ b=b_{(k)}dx^{(k)}.
$$
The \emph{(complete) contraction multiplication} is defined by
\begin{equation}\label{A.abG}
\Ga(T_\tau^\sa M)\times\Ga(T_\rho^0 M)
\ra\Ga(T_\tau^{\sa-\rho} M),
\qquad (a,b)\mt a\btdot b
\end{equation}
with
\begin{equation}\label{A.ab}
a\btdot b:=a_{(j)}^{(\ell)(k)}b_{(k)}\frac\pl{\pl x^{(\ell)}}
\otimes dx^{(j)}
\end{equation}
(so that $(k)$~runs through~$\BJ_\rho$). We infer from
\eqref{A.ab} that \eqref{A.abG} is a bilinear map
$$
C^\iy(T_\tau^\sa M)\times C^\iy(T_\rho^0M)
\ra C^\iy(T_\tau^{\sa-\rho}M).
$$

\smallskip
Let $g$ be a Riemannian metric on~$TM$. We write
\hb{g_\flat\colon TM\ra T^*M} for the (fiber-wise) Riesz
isomorphism. Thus
\hb{\dl g_\flat X,Y\dr=g(X,Y)} for $X$ and~$Y$ in $\Ga(TM)$ where
\hb{\pw\colon\Ga(T^*M)\times\Ga(TM)\ra\Ga(M)}
is the natural (fiber-wise) duality pairing. Hence
\hb{T_p^*M=(T_pM)'} for
\hb{p\in M}. The inverse
of~$g_\flat$ is
denoted by~$g^\sh$. Then~$g^*$, the adjoint Riemannian metric on~$T^*M$,
is defined by
\hb{g^*(\al,\ba):=g(g^\sh\al,g^\sh\ba)} for
\hb{\al,\ba\in\Ga(T^*M)}. In local coordinates
$$
g=g_{ij}\,dx^i\otimes dx^j,
\quad g^*=g^{ij}\,\frac\pl{\pl x^i}\otimes\frac\pl{\pl x^j},
$$
$[g^{ij}]$~being the inverse of the
\hb{(m\times m)}-matrix~$[g_{ij}]$. We also employ the notation
\hb{\prsn_g:=g}.

\smallskip
The metric~$g$ induces a vector bundle metric on~$T_\tau^\sa M$ which
we denote by~$g_\sa^\tau$. In local coordinates
$$
g_\sa^\tau(a,b)=g_{(i)(j)}g^{(k)(\ell)}a_{(k)}^{(i)}b_{(\ell)}^{(j)},
\ a,b\in\Ga(T_\tau^\sa M),
$$
where
$$
g_{(i)(j)}:=g_{i_1j_1}\cdots g_{i_\sa j_\sa},
\quad g^{(k)(\ell)}:=g^{k_1\ell_1}\cdots g^{k_\tau\ell_\tau}
$$
for
\hb{(i),(j)\in\BJ_\sa} and
\hb{(k),(\ell)\in\BJ_\tau}. Note
\hb{g_1^0=g} and
\hb{g_0^1=g^*} and
\hb{g_0^0(a,b)=ab} for
\hb{a,b\in\Ga(M)}. Moreover,
$$
\vsdot_{g_\sa^\tau}\colon\Ga(T_\tau^\sa M)\ra(\BR_+)^M,
\quad a\mt\sqrt{g_\sa^\tau(a,a)}
$$
is the \emph{vector bundle norm} on~$T_\tau^\sa M$ induced by~$g$.
It follows that the complete contraction is a fiber-wise continuous
bilinear map.

\smallskip
The Euclidean metric on~$\BR^m$ is named
\hb{g_m=(dx^1)^2+\cdots+(dx^m)^2}, and
\hb{\vsdot_m:=\vsdot_{g_m}} is the Euclidean norm.

\smallskip
The Levi--Civita connection on~$TM$ is denoted by
\hb{\na=\na_g}. We use the same symbol for its natural extension
to a metric connection on~$T_\tau^\sa M$. Then the corresponding covariant
derivative is the linear map
$$
\na\colon C^\iy(T_\tau^\sa M)\ra C^\iy(T_{\tau+1}^\sa M),
\quad a\mt\na a,
$$
defined by
\hb{\dl\na a,b\otimes X\dr:=\dl\na_Xa,b\dr} for
\hb{b\in C^\iy(T_\sa^\tau M)} and
\hb{X\in C^\iy(TM)}. It is a well-defined continuous linear
map from~$C^1(T_\tau^\sa M)$ into~$C(T_{\tau+1}^\sa M)$, as follows from
its local representation with the Christoffel symbols. For
\hb{k\in\BN} we define
$$
\na^k\colon C^\iy(T_\tau^\sa M)\ra C^\iy(T_{\tau+k}^\sa M),
\quad a\mt\na^ka
$$
by
\hb{\na^0a:=a} if
\hb{\sa=\tau=0}, and
\hb{\na^{k+1}:=\na\circ\na^k}, where
\hb{\na^1=d}, the differential on~$C^\iy(M)$. If
\hb{\Mg=(\BR^m,g_m)}, then
\hb{\na_{g_m}=\pl=(\pl_1,\ldots,\pl_m)},
the (Fr\'echet) derivative on~$\BR^m$.

\smallskip
Let $M_1$ and~$M_2$ be \hbox{$m$-dimensional} manifolds and
\hb{f\colon M_1\ra M_2} a diffeomorphism, in symbols:
\hb{f\in\Diff(M_1,M_2)}. The \emph{push-forward} of functions,
\hb{f_*\colon C^\iy(M_1)\ra C^\iy(M_2)}, is the linear bijection defined by
\hb{f_*u:=u\circ f^{-1}} for
\hb{u\in C^\iy(M_1)}. Its inverse, the \emph{pull-back}~$f^*$, is given by
\hb{f^*v=v\circ f} for
\hb{v\in C^\iy(M_2)}. The push-forward (by~$f$) of vector fields
is the vector space isomorphism
$$
f_*\colon C^\iy(TM_1)\ra C^\iy(TM_2),
\quad X\mt f_*X
$$
specified by
$$
(f_*X)(q)
:=(T_{f^{-1}(q)}f)X\bigl(f^{-1}(q)\bigr),
\  q\in M_2.
$$
Here
\hb{T_pf\colon T_pM_1\ra T_{f(p)}M_2} is the tangent map of~$f$ at
\hb{p\in M_1}. The pull-back of vector fields (by~$f$),
$$
f^*\colon C^\iy(TM_2)\ra C^\iy(TM_1),
\ Y\mt f^*Y,
$$
is given by
$$
(f^*Y)(p):=(T_pf)^{-1}Y\bigl(f(p)\bigr),
\ p\in M_1.
$$

\smallskip
The vector space isomorphisms
$$
f_*\colon C^\iy(T^*M_1)\ra C^\iy(T^*M_2),
\quad \al\mt f_*\al,
$$
introduced by
$$
(f_*\al)(q):=(T_qf^{-1})'\al\bigl(f^{-1}(q)\bigr),
\ q\in M_2,
$$
and
$$
f^*\colon C^\iy(T^*M_2)\ra C^\iy(T^*M_1),
\ \ba\mt f^*\ba,
$$
with
$$
(f^*\ba)(p):=(T_pf)'\ba\bigl(f(p)\bigr),
\ p\in M_1,
$$
are the push-forward and the pull-back, respectively, of covector fields.

\smallskip
Now we define
$$
f_*\colon C^\iy(T_\tau^\sa M_1)\ra C^\iy(T_\tau^\sa M_2)
$$
inductively by
$$
f_*(a\otimes b):=f_*a\otimes f_*b,
\ a\in C^\iy(T_{\tau-j}^{\sa-i}M),
\ b\in C^\iy(T_j^iM),
$$
where
\hb{i,j\in\{0,1\}} with
\hb{\sa\geq i} and
\hb{\tau\geq j}. An analogous definition applies to
$$
f^*\colon C^\iy(T_\tau^\sa M_2)\ra C^\iy(T_\tau^\sa M_1).
$$
Then $f_*$~is a vector space isomorphism and
\hb{(f_*)^{-1}=f^*}. Furthermore, if
\hb{f_1\in\Diff(M_1,M_2)}, then
\hb{(f_1\circ f_2)_*=f_{1\,*}\circ f_{2\,*}}.

\smallskip
Suppose that $(M_1,g_1)$ and~$(M_2,g_2)$ are Riemannian manifolds. Then
$$
f\colon(M_1,g_1)\ra(M_2,g_2)
$$
is an \emph{isometric diffeomorphism} if
\hb{f\in\Diff(M_1,M_2)} and
\hb{f_*g_1=g_2}.

\smallskip
We write~$dV_g$ for the Riemann--Lebesgue volume measure on~$\Mg$.
In local coordinates
\hb{\ka=(x^1,\ldots,x^m)} it is represented by
\hb{\ka_*dV_g=(\ka_*\sqrt{g})dx}, where
\hb{\sqrt{g}:=\bigl(\det[g_{ij}]\bigr)\vph{)}^{1/2}} and $dx$~is the Lebesgue
measure on~$\BR^m$.

\smallskip
In Section~\ref{sec-M} there occur product manifolds of the type
\hb{Z:=(0,1]\times B}, where $B$~is an
\hb{(m-1)}-dimensional manifold. If
\hb{\pl B\neq\es}, then
\hb{\{1\}\times\pl B} is~a `corner' of~$Z$, which is locally diffeomorphic to
\hb{[0,1)^2\times(-1,1)^{m-2}}. Everything said above has straightforward
 extensions to such manifolds with corners. Detailed investigations
 of manifolds with corners are given in~\cite{Ama25a}.

\end{document}